\documentclass[review,10pt]{elsarticle}
\usepackage[numbers]{natbib}
 \usepackage{mathrsfs}
\usepackage{graphicx}
\usepackage{epstopdf}
\usepackage{amssymb}
\usepackage{color}
\usepackage{amsmath}
\usepackage{float}
\usepackage{multirow}
\makeatletter\@addtoreset{equation}{section}

\newtheorem{thm}{Theorem}[section]
\newtheorem{cor}[thm]{Corollary}
\newtheorem{lem}[thm]{Lemma}

\newtheorem{rem}[thm]{Remark}

\newtheorem{algo}[thm]{Algorithm}

\numberwithin{equation}{section}
\numberwithin{table}{section}
\numberwithin{figure}{section}

\makeatletter\@addtoreset{equation}{section}

\newenvironment{proofed}[1]{\par \textbf{Proof}\quad #1}{\hfill \textbf{} $\Box$ }

\newtheorem{lemma}{\bf{Lemma}}[section]
\newtheorem{example}{\bf{Example}}[section]

\numberwithin{equation}{section}
\numberwithin{table}{section}
\numberwithin{figure}{section}

\begin{document}

%\title{  Finite difference schemes for the space-fractional diffusion equation with nonlinear time delay  \footnote{The research is
%supported by National Natural Science Foundation of China (No.
%11271068)}}

\begin{frontmatter}
	\title{An improved algorithm based on finite difference schemes for fractional boundary value problems with non-smooth solution}
	
	\tnotetext[label1]{This work was partially supported by the NSF of China
		( No. 11271068). The first author was also partially supported by the
		National University Student Innovation Program (No.1410286047) and  China Scholarship Council (No. 201506090065). }
	\author[seu]{Zhao-Peng Hao}
	\ead{zhaopenghao2015@gmail.com}
	\author[seu]{Wan-Rong Cao\corref{cor1}}%\fnref{label2}}
	\ead{wrcao@seu.edu.cn}
	\cortext[cor1]{Corresponding author. Tel/Fax: +86 25 52090590}

	%\fntext[label2]{}
	%\author[seu]{Zhizhong Sun}
	%\ead{zzsun@seu.edu.cn}
	
	\address[seu]{Department of Mathematics, Southeast University, Nanjing 210096, P.R.China.}
	
	\begin{abstract}
	In this paper, an efficient algorithm is presented by the extrapolation technique to improve the accuracy of finite difference schemes  for solving the fractional boundary value problems with non-smooth solution. Two popular finite difference schemes,  the weighted shifted Gr\"{u}nwald difference (WSGD) scheme and the fractional centered difference (FCD) scheme, are revisited and the error estimate of the schemes  is  provided in maximum norm.  Based on the analysis of leading singularity of exact solution for the underlying problem, it is demonstrated that,  with the use of the proposed algorithm, the improved WSGD  and FCD schemes can recover the second-order accuracy for non-smooth solution. Several numerical examples are given to validate our theoretical prediction. It is shown that both accuracy and convergence rate of numerical solutions can be significantly improved  by using the proposed algorithm.
%, particularly for the fractional boundary value problems with one-sided fractional derivative or two-sided symmetric derivatives and the time-dependent problems.

	\end{abstract}

\end{frontmatter}

\section{Introduction}
The aim of this work is to present an  efficient  numerical approach, which is based on the finite difference method, to solve space fractional diffusion equations (SpFDEs) with  non-smooth solutions.

In recent decades, anomalous diffusion has been widely considered in the investigation of transport dynamics in complex systems, such as underground environmental problem \cite{Hatano98}, fluid flow in porous materials \cite{BensonSMW2001},  anomalous transport in biology \cite{HoflingFran2013}, etc. SpFDEs can provide an adequate and accurate description of the super-diffusion process \cite{Meers2000,Metzler2000}.

With the increasing application of SpFDE in modeling problems which exhibit super-diffusion, how to obtain its accurate numerical solution has attracted considerable attention. A great number of numerical methods have been developed in literature, among which the finite difference method is one of the most popular and powerful methods. A shifted Gr\"{u}nwald formula for the Riemann-Liouville fractional derivative  was firstly proposed in \cite{MeerTad2004} to solve the space fractional differential equations, which is of first-order accuracy and leads to unconditionally stable schemes. Based on this work, some high-order finite difference schemes for SpFDE have been proposed subsequently, e.g.,  the  second-order extrapolation method \cite{TadjeranMS2006},   a class of weighted shifted  Gr\"{u}nwald formulas  \cite{HaoSC2015,TianZD2015}. As another popular scheme,  the  second-order fractional centered difference scheme  was presented in \cite{Celik-Duman2012} to solve the equation with Riesz fractional derivatives.
In all aforementioned works, however,  the convergence rates are obtained under the requirement of high regularity of the solution. Though the assumption is natural for canonical  partial differential equations with integer-order derivatives, it is too idealistic to satisfy for fractional differential equations (FDEs) in application.  In fact, the fractional derivatives are defined with weak singular kernel and the solution of FDEs inherits the weak singularity. Even smooth data cannot ensure smoothness of the solution \cite{Diethelm2010,HaoPLC2016,JinLPR2015, MaoCS2016}. In addition, the above-mentioned high-order schemes require that the solution and its first or even up to  higher-order derivatives have vanishing values at the boundary. When solving FDEs whose solutions  have neither high regularity near the boundary nor vanishing derivatives at the boundary, the schemes based on these ideal assumptions will actually  lead to numerical solutions of very low accuracy.

As to time-fractional initial value  problems, for which the weak singularity of solutions usually exists at the origin, several approaches have been proposed to deal with the weak singularity in order to obtain numerical solutions of uniformly high-order accuracy, such as adopting adaptive grids (nonuniform grids) to keep errors small near the  initial time \cite{Mclean2007,YusteJoaquin12,ZhangSunLiao2014},  or employing non-polynomial basis functions to include the correct singularity index \cite{CaoHX2003,FordMR2013}, or adding the correction terms to remedy the loss of accuracy and  recover high-order schemes \cite{CaoZZK2016,Lub86,Zeng2015, ZengZK2015}.  As for space fractional boundary value problems/or initial boundary value problems,
the solutions  generally have weak singularity near the boundary or  the end-points of both sides in one-dimensional case. So far, to our knowledge, only a few works have been presented to  numerical methods for SpFDEs with non-smooth solutions. Zhao and Deng \cite{Zhao-Deng2016} derived finite difference schemes on non-uniform meshes to increase the accuracy for solving SpFDEs with non-smooth solution. Jin and Zhou \cite{JinZ2015}  proposed a  singularity reconstruction strategy   to enhance
the convergence rate and gave a new finite element method for approximating
boundary value problems with Riemann-Liouville fractional derivatives.
 Mao and Shen \cite{MaoCS2016} developed a spectral Petrov-Galerkin method for FDEs with Riesz fractional derivatives, in which the error estimate in non-uniformly weighted Sobolev space shows that the errors decay exponentially even though the solution has singularities at the endpoints.

The main contribution of this work is to present an improved algorithm based on finite difference methods, which is readily implemented, applicable to  various kinds of finite difference schemes, and able to significantly enhance the accuracy of numerical solutions for SpFDEs with weak singularities.
To deal with the singularity and obtain  second-order accuracy, we first separate the solution $u$ of the considered problem into a regular/smooth part $u^r$ and a  singular/non-smooth part $\xi^s u^s$, where $\xi^s$ is the coefficient of the singular part.
 Then we adopt extrapolation   and posterior error correction techniques to approximate $\xi^s$ and recover the second-order accuracy of numerical solutions.

 Compared to  finite difference schemes on non-uniform meshes \cite{Zhao-Deng2016}, the proposed algorithm holds the  Toeplitz-like structure of the finite difference schemes, which  is a remarkable feature to allow low storage and the use of fast algorithms; see \cite{Pang-Sun2012,Wangh2013}.
  Although the proposed algorithm will cause the extra cost,   the increase of storage and computational cost is acceptable.
Numerical examples show that the improved WSGD and FCD schemes, which are obtained by applying the proposed algorithm to the WSGD and FCD schemes respectively, can produce more accurate numerical solutions than the corresponding original schemes without using the proposed algorithm; see Examples \ref{exm1}-\ref{exm3}. Even if the regularity of solution is unknown, that is, the ``singular part" of solution is given by the basic analysis of FDEs and some conjectures, we can still obtain satisfactory accuracy; see Examples \ref{exm1}-\ref{exm2} (Case II).

The rest of the paper is organized as follows. In Section 2, we introduce some necessary definitions and notations. Moreover, the WSGD scheme and the FCD scheme for SpFDEs are presented in this section.  To derive the algorithm, we carry out the error estimate in maximum norm of the two schemes in Section 3. The main algorithm is derived in Section 4. In Section 5, we give some numerical examples for solving fractional boundary value problems and SpFDEs with non-smooth solutions to illustrate the efficiency of the proposed algorithm. Finally, we give some concluding remarks.

\section{The second-order finite difference schemes for SpFDEs}

The one-dimensional SpFDE has the form
\begin{subequations}
\begin{eqnarray}
&&u_t-\theta \ _aD_x^{\beta }u-(1-\theta)\ _xD_b^{\beta }u=g(x,t), \; x\in(a,b),\; 0<t\leq T, \label{eq:Time-dependent-1}\\
&&u(x,0)=u^0(x), \;x\in[a,b],\label{eq:Time-dependent-2}\\
&&u(a,t)=u(b,t)=0, \;0<t\leq T, \label{eq:Time-dependent-3}
\end{eqnarray}
\end{subequations}
where $g$ is a given function, $\beta \in(1,2)$, and the notations $_aD_x^{\beta }$ and $_xD_b^{\beta }$ refer to the left-sided and right-sided Riemann-Liouville derivatives of order $\beta $ defined in \eqref{def-RL-1}; $\theta\in[0,1]$ is a parameter. In cases of $\theta=1, 0, 1/2$,  \eqref{eq:Time-dependent-1} is known as the SpFDE with the left-sided Riemann-Liouville derivative, the right-sided Riemann-Liouville derivative and the Riesz fractional derivative, respectively.

\subsection{Preliminaries}

We first briefly recall the Riemann-Liouville fractional integrals and derivatives.

The  left-sided and right-sided   Riemann-Liouville fractional integrals of the function $v(x)$ are respectively  defined by
\begin{equation}\label{def-1}
_aD_x^{-\beta } v(x) =\frac{1}{\Gamma(\beta )}\int_a^x \frac{v(\zeta) }{(x-\zeta)^{1-\beta }} d\zeta ,\quad x>a, \quad \beta  \in (0,1),
\end{equation}
\begin{equation}\label{def-1R}
_xD_b^{-\beta } v(x) =\frac{1}{\Gamma(\beta )}\int_x^b \frac{v(\zeta) }{(\zeta-x)^{1-\beta }} d\zeta ,\quad x<b, \quad \beta  \in (0,1).
\end{equation}
\iffalse
It is not hard to check that the formulas
\begin{equation}\label{integral-formula}
_aD_x^{-\sigma}(x-a)^{\xi}= \frac{\Gamma(\xi+1 )}{\Gamma(\xi+1+\sigma)}(x-a)^{\xi+\sigma},\quad
\ _xD_b^{-\sigma}(b-x)^{\xi}= \frac{\Gamma(\xi+1 )}{\Gamma(\xi+1+\sigma)}(b-x)^{\xi+\sigma}
\end{equation}
hold for any $x\in (a,b) ,$ $\sigma\in (0,1) $ and $\xi >-1.$
The formulas in \eqref{integral-formula} will be used to calculate the functions explicitly in the numerical tests.
\fi

The $\beta  ~(n-1< \beta  < n)$ order left-sided and right-sided  Riemann-Liouville fractional derivatives of the function $v(x)$ on $[a,b]$ are defined by
\begin{equation}\label{def-RL-1}
_aD_x^{\beta }v(x)=D^n\ _aD_x^{\beta -n}  v(x),\quad _xD_b^{\beta }v(x)=(-1)^nD^n  \ _xD_b^{\beta -n} v(x).
\end{equation}
where $D^{n}:=d^{n}/dx^n.$
If $\beta =n,$ then $_0D_x^{\beta }v(x)=D^{n}v(x)$ and $_xD_b^{\beta }v(x)=(-1)^nD^{n}v(x).$

\iffalse
The Riesz potential i.e, Riesz fractional integral for  $\beta  \in (0,1)$ is defined as

\begin{eqnarray}
\ _a^1\mathbb{D}_b^{-\beta } v(x) &=&\frac{1}{2\cos(\beta  \pi/2)\Gamma(\beta )} \int_a^b \frac{sign(x-t) v(\xi) }{|x-\xi|^{1-\beta }} d\xi \nonumber\\
&=&\frac{1}{2\cos(\beta  \pi/2)}[\ _aD_x^{-\beta } v(x)-\ _xD_b^{-\beta } v(x)] ,\quad a<x<b,\label{def-Reisz-1}\\
\ _a^2\mathbb{D}_b^{-\beta } v(x) &=&\frac{1}{2\cos(\beta  \pi/2)\Gamma(\beta )} \int_a^b \frac{v(\xi) }{|x-\xi|^{1-\beta }} d\xi \nonumber\\
&=&\frac{1}{2\cos(\beta  \pi/2)}[\ _aD_x^{-\beta } v(x)+\ _xD_b^{-\beta } v(x)],\quad a<x<b . \label{def-Reisz-2}
\end{eqnarray}
where $sign$ is sign function
For  $\beta  ~(n-1< \beta  < n),$  Riesz fractional derivative of order $\beta $ is defined as
\begin{equation}\label{def-Reisz-derivative}
_a\mathbb{D}_b^{\beta } v(x):=
\left\{  \begin{array}{l}
D^n\ _a^1\mathbb{D}_b^{\beta -n} v(x), \quad n ~is~ odd,\\
D^n\ _a^2\mathbb{D}_b^{\beta -n} v(x), \quad n ~is ~even,\\
\end{array}
\right.
\quad a<x<b.
\end{equation}
\fi
If we take $\theta=1/2$ in  \eqref{eq:Time-dependent-1}, then the two-sided fractional derivative  will be symmetric and closely related to the Riesz fractional derivative, which is defined by
\begin{equation}\label{eq:riesz}
_a\mathbb{D}_b^{\beta } v(x)= -\frac{1}{2\cos( \frac{\beta  \pi}{2})}[\ _aD_x^{\beta }v(x)+\ _xD_b^{\beta }v(x)],\;1<\beta <2.
\end{equation}
It should be noted that $a$ (or $b$ ) in \eqref{def-1}-\eqref{eq:riesz} is allowed to take $-\infty$ (or $+\infty$), so that the fractional integrals and derivatives can be defined on the whole real axis. For the details of definitions and properties of the fractional integrals and derivatives, we refer to \cite{Samko1993}.
\iffalse
In the end of this subsection, we provide  a useful formula, which will be required to derive the improved algorithm in Section 4.
\begin{equation}\label{eq-theta-3}
_a\mathbb{D}_b^{\beta } [(x-a)^{\frac{\beta }{2}}(b-x)^{\frac{\beta }{2}}]=-\Gamma(\beta +1),\quad  1<\beta  \leq 2.
\end{equation}
\fi

%As another commonly used fractional derivative,  Riesz fractional derivative with order $\beta $ of function $v(x)$ is defined as
% \begin{equation}\label{fs2}
% \mathbb{D}^{\beta }v(x)=-\frac{1}{2\cos(\beta  \pi/2)\Gamma(2-\beta )} \frac{d^2}{dx^2}\int_{a}^{b}\frac{v(\xi)}{|x-\xi|^{\beta -1}}d\xi,\quad 1<\beta < 2,
% \end{equation}

\subsection{Second-order approximations of fractional derivatives}
In this part, we introduce the second-order weighted shifted Gr\"{u}nwald difference approximation \cite{HaoSC2015} for the Riemmann-Liouville fractional derivatives and the fractional centered difference approximation \cite{Celik-Duman2012} for the Riesz fractional derivative.

\begin{lem}(See \cite{HaoSC2015})\label{lem1}
	Let  $v\in L_1{(\mathbb{R})},$  ${_{-\infty}}D_x^{\beta +2}v(x) ,$ $_xD_{+\infty}^{\beta +2}v(x) $ and their Fourier transforms belong to $L_1{(\mathbb{R})}.$  Then for the given step size $h$,  the left- and right-sided Riemann-Liouville fractional derivatives with order $\beta $ of $v$ at each point $x$ can be approximated by the  weighted shifted Gr\"{u}nwald difference approximation with second-order accuracy,
	\begin{eqnarray}
	{_{-\infty}}D_x^{\beta }v(x)&=&\frac{1}{h^{\beta }} \sum_{k=0}^{+\infty} w_k^{(\beta )}v(x-(k-1)h)+\mathcal{O}(h^2),\label{a17}\\
	_xD_{+\infty}^{\beta }v(x)&=&\frac{1}{h^{\beta }} \sum_{k=0}^{+\infty} w_k^{(\beta )}v(x+(k-1)h)+\mathcal{O}(h^2),\label{a18}
	\end{eqnarray}
where
\begin{equation}\label{a21}
w_0^{(\beta )}=\lambda_1g_0^{(\beta )},\; w_1^{(\beta )}=\lambda_1g_1^{(\beta )}+\lambda_0g_0^{(\beta )},\;
w_k^{(\beta )}=\lambda_1g_k^{(\beta )}+\lambda_0g_{k-1}^{(\beta )}+\lambda_{-1}g_{k-2}^{(\beta )},\; k\geq 2,
\end{equation}
and the weights $\{g_k^{(\beta )}\}$ in \eqref{a21} are the coefficients of the power series of function $(1-z)^{\beta },$ i.e.,
\begin{equation}\label{a3}
(1-z)^{\beta }=\sum_{k=0}^{+\infty}(-1)^k\binom{\beta }{k}z^k=\sum_{k=0}^{+\infty} g^{(\beta )}_{k}z^k,\;-1<z\leq 1,
\end{equation}
and $$\lambda_1=\frac{\beta ^2+3\beta +2}{12},\quad \lambda_0=\frac{4-\beta ^2}{6},\quad \lambda_{-1}=\frac{\beta ^2-3\beta +2}{12}.$$
\end{lem}
Note that $\{g_k^{(\beta)}\}$ can be computed recursively, that is,
$$g_{k+1}^{(\beta)}=(1-\frac{\beta+1}{k+1})g_k^{(\beta)},\quad k=0, 1,  2,\ldots.. $$
\begin{lem}(See \cite{Celik-Duman2012})\label{lem2}
	Denote
	$$\widetilde{w}^{(\beta )}_k= -\frac{(-1)^k\Gamma(\beta +1)}{\Gamma(\frac{\beta }{2}-k+1)\Gamma(\frac{\beta }{2}+k+1)},\;k=\pm 1, \pm 2,\ldots.$$
	Let  $v\in L_1{(\mathbb{R})},$  ${_{-\infty}}\mathbb{D}_{+\infty}^{\beta +2}v(x) $ and its Fourier transform belong to $L_1{(\mathbb{R})}.$ Then for a given step size $h$,  it holds that
	%\begin{equation}\label{fgl5}
	%\frac{1}{h^{\beta }}\sum_{k=-\infty}^{+\infty}\widetilde{g}^{(\beta )}_kv(x-kh)=\mathbb{D}^{\beta }v(x)+\sum_{l=1}^{n/2-1}\biggl[c_l^{\beta }\cdot(-1)^l\mathbb{D}^{\beta +2l}v(x)\biggr]h^{2l}+O(h^{n})
	%\end{equation}
	\begin{equation}\label{fgl5}
	\ _{-\infty}\mathbb{D}_{+\infty}^{\beta } v(x)=\frac{1}{h^{\beta }}\sum_{k=-\infty}^{+\infty}\widetilde{w}^{(\beta )}_kv(x-kh)
	+\mathcal{O}(h^{2})
	\end{equation}
	uniformly for  $x\in \mathbb{R}.$
\end{lem}

The weights $\{\widetilde{w}_k^{(\beta )}\}$ in above lemma are the coefficients of Fourier series  of the function $ |2\sin(\frac{z}{2})|^{\beta },$ i.e.,
\begin{equation*}
|2\sin(\frac{z}{2})|^{\beta }=\sum_{k=-\infty}^{+\infty}\widetilde{w}_k^{(\beta )}e^{ikz}.
\end{equation*}
Noting $\Gamma(z+1)=z\Gamma(z),$ we can also  write $\widetilde{w}^{(\beta )}_k$ in the recursive way
\begin{equation*} \widetilde{w}^{(\beta )}_k=(1-\frac{\beta +1}{\frac{\beta }{2}+k} )\widetilde{w}^{(\beta )}_{k-1},\quad k=\pm 1, \pm 2,\ldots.  \end{equation*}
%where
%$\widetilde{w}^{(\beta )}_0= -\frac{\Gamma(\beta +1)}{\Gamma(\frac{\beta }{2}+1)\Gamma(\frac{\beta }{2}+1)}.$

\subsection{Derivation of the difference scheme}

Take an integer $M.$ Let $I_h\equiv \{ x_j~|~0\leq j\leq M\}$ be a uniform mesh of the interval $[a,b],$ where $x_j=a+jh,~0\leq j\leq M$ with $h=(b-a)/M.$ Suppose  $v=\{v_{j} \}$ is a grid function on $I_h.$   We define the left- and right-sided weighted shifted Gr\"{u}nwald finite difference operators
\begin{eqnarray}\label{eq:wsgd}
&&  \delta_{x,-}^{\beta }v_j=\frac{1}{h^{\beta }}\sum_{k=0}^{j} w_k^{\beta }v_{j-k+1}, \quad \delta_{x,+}^{\beta }v_j=\frac{1}{h^{\beta }}\sum_{k=0}^{M-j} w_k^{\beta }v_{j+k-1},\quad 1\leq  j \leq M-1,
\end{eqnarray}
and the fractional centered difference operator
\begin{eqnarray}\label{eq:centralapp}
&&  \Delta_{x}^{\beta }v_j=\frac{1}{h^{\beta }}\sum_{k=-M+j}^{j}  \widetilde{w}_{k}^{\beta }v_{j-k},\;1\leq  j \leq M-1 .
\end{eqnarray}

After a suitable time discretization of \eqref{eq:Time-dependent-1},
we are led to solve, at each time step, a fractional elliptic problem of the following kind:
\begin{subequations}
\begin{eqnarray}
&&  \alpha u-\theta\ _aD_x^{\beta }u-(1-\theta)\ _xD_b^{\beta }u=f(x), \quad x \in(a,b),\label{vp-eq-1} \\
&& u(a)=u(b)=0,\label{vp-b-1}
\end{eqnarray}
\end{subequations}
where $ \alpha$ is a positive scaling constant.

Define
\begin{equation*}
\widetilde{u}(x)= \left \{\begin{array}{cc}
u(x) ,& x\in[a,b],\\
0,~&  \mbox{otherwise}.
\end{array}
\right.
\end{equation*}
To derive the difference scheme for \eqref{vp-eq-1}-\eqref{vp-b-1}, we assume that the zero-extended solution $\widetilde{u}(x) $  satisfies the conditions in Lemma \ref{lem1}.

Consider \eqref{vp-eq-1} at the grid points $x=x_j$, we have
\begin{eqnarray}
&& \alpha u(x_j)-\theta \ _a D_x^{\beta }u(x_j) -(1-\theta) \ _x D_b^{\beta }u(x_j)=f(x_j).\label{Discrete-eq-1}
\end{eqnarray}
By Lemma \ref{lem1} and the definition \eqref{eq:wsgd}, we have
\begin{eqnarray}
&&\alpha u(x_j) -\theta\delta_{x,-}^{\beta }u(x_j)-(1-\theta)\delta_{x,+}^{\beta }u(x_j)=f(x_j)+R_j, \quad  1\leq j\leq M-1,
\label{Discrete-eq-2}
\end{eqnarray}
where there exists a constant $c_R$ independent of step size $h$ such that
\begin{eqnarray}\label{trunction-err}
|R_j|\leq  c_Rh^2.
\end{eqnarray}
Notice the homogenous Dirichlet boundary conditions \eqref{vp-b-1}.
Omitting $R_j$ in \eqref{Discrete-eq-2} and
denoting by $u_{j}$ the numerical approximation of $u(x_j),$  $f_j=f(x_j)$, we get the weighted shifted Gr\"{u}nwald finite difference scheme (abbreviated as \textbf{WSGD})
\begin{subequations}
\begin{eqnarray}
&& \alpha u_j-\theta \delta_{x,-}^{\beta }u_j- (1-\theta)\delta_{x,+}^{\beta }u_j=f_j, \quad 1\leq j\leq M-1,\label{scheme-1}\\
&& u_0=u_M=0.\label{scheme-2}
\end{eqnarray}
\end{subequations}
\iffalse
Denote
\begin{equation}\label{Stiffness-maxtrix}
S_{\beta }=\frac{1}{h^{\beta }}
\left(
\begin{array}{llllll}
w_1^{(\beta )} &w_0^{(\beta )}&~  & ~ &~&~  \\
w_2^{(\beta )} &w_1^{(\beta )}& w_0^{(\beta )} & ~ & ~&~ \\
w_3^{(\beta )} & w_2^{(\beta )} & w_1^{(\beta )} & \ddots & ~&~ \\
\vdots   &  \vdots  & \vdots  & \ddots &\ddots&~\\
w_{M-2}^{(\beta )} & w_{M-3}^{(\beta )} &  w_{M-4}^{(\beta )} &  \cdots &  w_{1}^{(\beta )} &  w_{0}^{(\beta )} \\
w_{M-1}^{(\beta )} &  w_{M-2}^{(\beta )} &  w_{M-3}^{(\beta )}& \cdots& w_{2}^{(\beta )} &  w_{1}^{(\beta )}\\
\end{array}
\right).
\end{equation}

To calculate the difference equations, we can recast the difference scheme \eqref{scheme-1}-\eqref{scheme-2}  into matrix  formulation,
$AU=F,$ with   $U=(u_j)_{1\leq j\leq M-1},$  $F=(f_j)_{1\leq j\leq M-1} $ and  $A=\theta S_{\beta }+(1-\theta) S_{\beta }^{T}+I,$ where $S_{\beta }^{T}$ is the transpose of matrix  $S_{\beta }$ and  $I$ denotes a unit matrix.
\fi
Particularly,  by \eqref{eq:riesz}, when $\theta=\frac{1}{2},$  \eqref{vp-eq-1}  reduces to the following equation with the Riesz fractional derivative
\begin{subequations}\label{vp-R}
\begin{eqnarray}
&&\alpha u(x)+\cos(\frac{\beta  \pi}{2}) _a\mathbb{D}_b^{\beta } u(x)=f(x), \quad x \in(a,b),\label{vp-R-1} \\
&& u(a)=u(b)=0.\label{vp-R-2}
\end{eqnarray}
\end{subequations}
For \eqref{vp-R-1}-\eqref{vp-R-2}, besides using the scheme \eqref{scheme-1}-\eqref{scheme-2}, we can also adopt the fractional centered approximation  to discretize the Riesz fractional derivative directly, and obtain the finite centered difference scheme (abbreviated as \textbf{FCD})
\begin{subequations}
\begin{eqnarray}
&&\alpha u_j+ \cos(\frac{\beta  \pi}{2})\Delta_{x}^{\beta }u_j=f_j, \quad 1\leq j\leq M-1,\label{scheme-c1}\\
&& u_0=u_M=0.\label{scheme-c2}
\end{eqnarray}
\end{subequations}
\section{Analysis of the finite difference schemes}
We analyze the convergence and stability of the WSGD scheme \eqref{scheme-1}-\eqref{scheme-2} and the FCD scheme \eqref{scheme-c1}-\eqref{scheme-c2} in this part.

First, we need to introduce some necessary notations.
The set of infinite grid is denoted by $h\mathbb{Z},$ with grid points $x_j=jh$ for $j\in \mathbb{Z},$ the set of all integers.
For any grid functions  $u=\{u_j\}$, $v=\{v_j\}$ on $h\mathbb{Z} ,$ the discrete inner product and the associated norm are defined as
$$( u,v)=h\sum_{j\in \mathbb{Z}}u_j\bar{v}_j,\quad \|u\|^2=( u,u) , $$
and   the discrete maximum norm is denoted by
$$\|u\|_{\infty}=\sup_{j\in \mathbb{Z}} |u_j| .$$
Set $L_h^2:=\{u\,|\,u=\{u_j\},\,\|u \|< +\infty \}.$ For $u\in L_h^2,$ we define the semi-discrete Fourier transform \cite{Trefethen2000} $\hat{u}: [-\frac{\pi}{h}, \frac{\pi}{h} ]\rightarrow \mathbb{C}$ by
\begin{equation}
\hat{u}(k):= h \sum_{j\in \mathbb{Z}} u_j e^{-ikx_j}, \label{semift}
\end{equation}
and   the inverse  semi-discrete Fourier transform
\begin{equation}
u_j=\frac{1}{2\pi} \int_{-\pi/h}^{\pi/h} \hat{u}(k) e^{ikx_j}dk. \label{inverse-1}
\end{equation}
 It is not hard to check that  Parseval's identity
\begin{equation}\label{parseval}
(u,v)=\frac{1}{2\pi}\int_{-\pi/h}^{\pi/h} \hat{u}(k) \overline{\hat{v}(k)}dk
\end{equation}
holds.
For a fixed constant  $\sigma \in (0,1],$ we define the fractional Sobolev semi-norm $|\cdot|_{H^{\sigma}}$ and norm $\|\cdot\|_{H^{\sigma}}$ as
\begin{eqnarray}
&&|u|^2_{H^{\sigma}}=\int_{-\pi/h}^{\pi/h}|k|^{2\sigma}|\hat{u}(k)|^2dk,\label{sigma-1}\\ &&\|u\|^2_{H^{\sigma}}=\int_{-\pi/h}^{\pi/h}(1+|k|^{2\sigma})|\hat{u}(k)|^2dk.\nonumber
\end{eqnarray}
Obviously, $\|u\|_{H^{\sigma}}^2=\|u\|^2+|u|_{H^{\sigma}}^2. $ Set $H^{\sigma}_h:=\{u\,|\,u=\{u_j\},\,\|u \|_{H^{\sigma}}< +\infty \}.$
 Denote $\mathcal{V}_h=\{v\;|\;v=\{v_{j}\},~ 0\leq j\leq M\}$,  $ \overset{\circ}{\mathcal V}_h=\{v~|~v\in \mathcal V_h,  v_0=v_M=0 \}.$ It is readily to know that for any $v\in \overset{\circ}{\mathcal V}_h,$ the Parseval's identity \eqref{parseval} still holds. In fact, it is sufficient  to extend  $v$ to infinite sequence by setting  $v_j=0$ for $ j\neq 0,1\cdots, M. $

Before presenting the maximum-norm error estimate, we cite several necessary lemmas.
\begin{lemma}(See \cite{HaoFCZ2015,Hao-Sun2016})\label{Sobolev-1}
	For  $\frac{1}{2}< \sigma \leq 1,$ there exists a constant $c_{0}>0$  depending on the parameter $\sigma$  but  independent of $h>0$ such that
	$$\|u\|_{\infty} \leq c_{0}\|u\|_{H^{\sigma}} $$
	for all $u\in H_h^{\sigma}.$
\end{lemma}

The left- and right-sided  fractional difference operators are adjoint to each other, which states as follows.
\begin{lemma}(See \cite{HaoSC2015})\label{adjoint} For  $1< \beta  \leq 2,$ we have
	$$ (\delta_{x,+}^{\beta }u,u)=(u,\delta_{x,-}^{\beta }u). $$
\end{lemma}

In the next lemma, we have the  norm equivalence  which is essential to the analysis of the scheme.
\begin{lemma}(See \cite{HaoFCZ2015,Hao-Sun2016})\label{Sobolev-2} For  $1< \beta  \leq 2,$ we have
	$$c_{\beta } |u|^2_{H^{\frac{\beta }{2}}} \leq -(\delta_x^{\beta }u,u) \leq |u|^2_{H^{\frac{\beta }{2}}},$$
	where $\delta_x^{\beta }u=-\frac{1}{2\cos(\frac{\beta  \pi}{2})}(\delta_{x,+}^{\beta }u +\delta_{x,-}^{\beta }u)$ and $c_{\beta }=\frac{2^{\beta }(1-\beta ^2)}{3\pi^{\beta }\cos(\frac{\beta  \pi}{2})} .$
\end{lemma}

From above two lemmas, it follows that
\begin{eqnarray}\label{ineq1}
c_{\ast}|u|^2_{H^{\frac{\beta }{2}}} \leq -(\delta_{x,+}^{\beta }u,u)= -(\delta_{x,-}^{\beta }u,u) \leq - \cos(\frac{\beta  \pi}{2}) |u|^2_{H^{\frac{\beta }{2}}},
\end{eqnarray}
where $ c_{\ast}=\frac{2^{\beta}(\beta ^2-1)}{3\pi^{\beta }}.$

%Let
%$$V_h=\{v|v=(v_0,\cdots,v_M ),~v_0=v_M=0\}$$ be a grid function space on $I_h$.
%For any   real  $u, v\in V_h,$ the aforementioned properties still hold. In fact, it is sufficient  to extend  $u$ to infinite sequence by setting  $u_j=0$ for $ j\neq 0,\cdots M. $

\begin{thm}\label{convergence}
	The WSGD scheme \eqref{scheme-1}-\eqref{scheme-2} is uniquely solvable and of second-order convergence in maximum norm. More precisely,  suppose the zero-extended function $\widetilde{u}(x)\in  L_1{(\mathbb{R})},$   ${_{-\infty}}D_x^{\beta +2}\widetilde{u}(x) $, $_xD_{+\infty}^{\beta +2}\widetilde{u}(x) $ and their Fourier transforms belong to $L_1{(\mathbb{R})}$ as well.
	Then  there exists a constant $c$ such that
	$$\max_{1\leq j\leq M-1}|u(x_j)-u_j|\leq c h^2 $$
	holds for all $\beta  \in(1, 2].$
\end{thm}
\begin{proofed}
	Since the unique solvability  has been given in  \cite{HaoSC2015}, here we just focus our attention on the convergence. Define the error grid functions as follows:
	\begin{displaymath}
	e_j=u(x_j)-u_j,\;0\leq j\leq M.
	\end{displaymath}
	
	Subtracting \eqref{scheme-1} from \eqref{Discrete-eq-2} leads to the error equation
	\begin{eqnarray}
	&&\alpha e_j -\theta \delta_{x,-}^{\beta }e_j- (1-\theta)\delta_{x,+}^{\beta }e_j=R_j, \quad \;1\leq j\leq M-1.,\label{eq-error}\\
	&& e_0=e_M=0.\label{eq-error-2}
	\end{eqnarray}
	Taking the discrete inner product of \eqref{eq-error} with $e$  on both sides gives
	\begin{eqnarray}
	&& \alpha(e,e)-\theta (\delta_{x,-}^{\beta }e, e)- (1-\theta)(\delta_{x,+}^{\beta }e, e)=(R,e).\label{eq-error-3}
	\end{eqnarray}
	For the left hand side of \eqref{eq-error-3}, by  Lemma \ref{Sobolev-1}, \eqref{ineq1} and \eqref{eq-error-2}, we have
	\begin{eqnarray}
	\quad \alpha (e,e)-\theta (\delta_{x,-}^{\beta }e, e)- (1-\theta)(\delta_{x,+}^{\beta }e, e)
	&\geq&  \alpha \|e\|^2+ c_{\ast}|e|^2_{H^{\frac{\beta }{2}}}
	\geq \frac{c_{\ast}}{c_0}\|e\|^2_{\infty} . \label{eq-error-4}
	\end{eqnarray}
	For the right hand side of \eqref{eq-error-3}, using the Cauchy-Schwarz inequality yields
	\begin{eqnarray}
	(R,e)\leq \|R\| \|e\|\leq (b-a)\|R\|_{\infty} \|e\|_{\infty}.\label{eq-error-5}
	\end{eqnarray}
	Substituting \eqref{eq-error-4} and \eqref{eq-error-5} into \eqref{eq-error-3} gives
	\begin{eqnarray*}
	\frac{c_{\ast}}{c_0}\|e\|^2_{\infty} \leq (b-a)\|R\|_{\infty} \|e\|_{\infty}.\label{eq-error-6}
	\end{eqnarray*}
	Consequently,   we have
	\begin{eqnarray*}
	\frac{c_{\ast}}{c_0}\|e\|_{\infty} \leq (b-a)\|R\|_{\infty} \leq c_R(b-a)h^2.\label{eq-error-6}
	\end{eqnarray*}
	This completes the proof.	
\end{proofed}

Similar to the proof of  Theorem \ref{convergence}, the stability of the WSGD scheme \eqref{scheme-1}-\eqref{scheme-2} can be  obtained straightforwardly.
\begin{cor}
	The WSGD scheme \eqref{scheme-1}-\eqref{scheme-2} is unconditionally stable to   the right hand term.
\end{cor}

For the difference scheme \eqref{vp-R-1}-\eqref{vp-R-2}, we have the following result.
\begin{thm}\label{convergence-R}
	The FCD scheme \eqref{vp-R-1}-\eqref{vp-R-2} is uniquely solvable, unconditionally stable and of second-order convergence in maximum norm.
\end{thm}
As the proof of this theorem is similar to the proof of Theorem \ref{convergence}, we omit it here.

%\begin{rem}
%
%
%Let $u$ be the exact solution of problem \eqref{c1}, and $u_h$ be the solution of difference scheme \eqref{c5}-\eqref{c17}.  Suppose the zero  extended function $\widetilde{u}(x)\in  \mathscr{C}^{4+\beta }(\mathbb{R}) .$
%Then   the following estimate
%$$\|u(x_j)-u_j\|_{\infty}\leq c h^4 $$
%holds for all $1< \beta  \leq 2.$
%
%
%\end{rem}
\section{Improved algorithm }
Due to the feature  of  weak singularity kernel in the definitions of fractional derivatives, solutions of fractional equations naturally  inherit the characteristic of weak singularity. In this section, we first justify the   leading  weak singularity in terms of \eqref{vp-eq-1}-\eqref{vp-b-1} with one-sided and two-sided symmetrical fractional derivatives.   Then  we use  extrapolation technique and  the  posterior   error correction method to recover the second-order accuracy of the WSGD scheme \eqref{scheme-1}-\eqref{scheme-2} and the FCD scheme \eqref{scheme-c1}-\eqref{scheme-c2}.

\subsection{The representation of the solution}\label{sec:41}

For $\theta=1$ or $\theta=0$ in \eqref{vp-eq-1},  regularity of the equation has been investigated in \cite{JinLPR2015} and also demonstrated in \cite{HaoPLC2016}, and we introduce the results as follows.

When $\theta=1,$ the equation \eqref{vp-eq-1} reduces to
\begin{equation}\label{eq:leftcase}
-\ _aD_x^{\beta }u=f-\alpha u.
\end{equation}
Let $\widetilde{f}=f-\alpha u.$ Then integrating on both sides of \eqref{eq:leftcase} twice reads
{\begin{equation}
-\ _aD_x^{\beta -2}u=\ _aD_x^{-2}\widetilde{f}+C_1(x-a)+C_2,
\end{equation}}
where $C_1$ and $C_2$ are coefficients to be determined. Taking $x\rightarrow a^{+}$ leads to $C_2=0$ in above equality. Since $u(a)=0,$  performing the fractional derivative operator $\ _0D_x^{2-\beta }$ on both sides gives
$$u=-\ _aD_x^{-\beta }\widetilde{f}- \frac{C_1}{\Gamma(\beta)}(x-a)^{\beta -1}. $$
Thus, it is readily to know that  the  leading weak singularity term of the solution is $ C(x-a)^{\beta -1}.$

When $\theta=0$ in \eqref{vp-eq-1}, that is, the equation only contains the  right-sided fractional derivative, similarly, one can derive that the leading weak singular term should be of the form $C(b-x)^{\beta -1} $.
\iffalse
For the two-sided case, we have the generalized Abel integration equation as  belows
{\color{red}\begin{equation}\label{eq-theta}
-\theta \ _aD_x^{\beta -2}u-(1-\theta) \ _xD_b^{\beta -2}u=-\ _aD_x^{-2}\widetilde{f}+C_1x+C_2.
\end{equation}}
\fi

Unlike the one-sided case,  where the leading weak singular term of solution can be readily determined by the transformation between fractional integrals and derivatives, the two-sided fractional derivatives will be far more complicated.  In fact,  based on the solution representation theory in \cite{Samko1993},  the inverse of the two-sided  fractional integrals involves the composition of the left-sided and right-sided integrals and derivatives, which makes the process of seeking the leading weak singularity  extremely difficult. Fortunately, for the symmetrical case, i.e., $\theta=1/2$ in \eqref{vp-eq-1},  we can get the classical   Carleman integral  equation with  Reisz potential, i.e.,
\begin{equation}\label{eq-theta-2}
-\frac{1}{2\Gamma(2-\beta)}\int_{a}^b \frac{u(\zeta)}{|x-\zeta|^{\beta -1}}d\zeta  =-\ _aD_x^{-2}\widetilde{f}+C_1x+C_2.
\end{equation}
By the spectral  relationship with the Reisz potential and the homogeneous boundary condition, one can find that the leading weak singular term has the form $ C (x-a)^{\frac{\beta }{2}}(b-x)^{\frac{\beta }{2}}.$ Since  the  derivation involves the theory of orthogonal Jacobi polynomials which will deviate from the main idea of this paper, herein we do not dwell on it. We refer to \cite{MaoCS2016} for the details. In next section, we will verify our prediction of singularities numerically; see Examples \ref{exm1}-\ref{exm2} (Case II).

\subsection{The posterior error correction method}
For ease of presentation, we introduce operators $D_{x,\theta}^{\beta }$ and $\delta_{x,\theta}$
$$D_{x,\theta}^{\beta }:=\theta \ _a D_x^{\beta }+(1-\theta)\ _x D_b^{\beta }, \quad \delta_{x,\theta}:=\theta\delta_{x,-}^{\beta }+(1-\theta)\delta_{x,+}^{\beta }.$$

Based on the discussion in subsection \ref{sec:41} and taking the homogeneous boundary conditions into account, it is reasonable to assume that the solution to \eqref{vp-eq-1} has the form
\begin{equation}\label{split-part}
u(x)=u^r(x)+\xi^s u^s(x),
\end{equation}
where  $u^r(x)$ is a regular part of $u(x)$ and its zero-extended function $\widetilde{u}^r(x)$ satisfies the assumptions in Lemma \ref{lem1}, $\xi^s$ is a scaling constant and $u^s$ is  of the form \begin{equation}\label{eq:singularterm}
u^s(x)=(x-a)^{\rho(\theta,\beta )}(b-x)^{\rho(1-\theta,\beta )},
\end{equation}
where $\rho$ is a given non-negative function related to $\theta$ and $\beta.$
%, and $$\displaystyle\max_{\theta\in[0,1],\,\beta\in(1,2)}\rho(\theta,\beta)<2.$$
 Then, there exists a  constant $\sigma<2$, which is dependent on $\rho$, but independent of $h$, such that
\begin{eqnarray}\label{eq:singularrate}
D_{x,\theta}^{\beta }u^s(x_j)= \delta_{x,\theta}u^s(x_j)+\mathcal{O}(h^{\sigma}).
\end{eqnarray}
For the smooth part $u^r(x)$, by Lemma \ref{lem1}, we have
\begin{eqnarray}\label{eq:smoothrate}
D_{x,\theta}^{\beta }u^r(x_j)= \delta_{x,\theta}u^r(x_j)+\mathcal{O}(h^2).
\end{eqnarray}
By \eqref{split-part} and \eqref{eq:smoothrate}, we have
\begin{eqnarray}
\quad D_{x,\theta}^{\beta }u(x_j)&=&D_{x,\theta}^{\beta }u^r(x_j)+\xi^s\cdot D_{x,\theta}^{\beta }u^s(x_j)\nonumber\\
&=&\delta_{x,\theta}u^r(x_j)+\xi^s\cdot D_{x,\theta}^{\beta }u^s(x_j)+\mathcal{O}(h^2)\nonumber\\
&=&\delta_{x,\theta}u(x_j)+\xi^s [D_{x,\theta}^{\beta }u^s(x_j)-\delta_{x,\theta}u^s(x_j)]+\mathcal{O}(h^2).
\end{eqnarray}

To simplify the description, we define an operator
$$\mathcal{L}u:= \alpha u-\theta\ _aD_x^{\beta }u-(1-\theta)\ _xD_b^{\beta }u$$
and let $u_h,~ u^r_h,~ u^s_h,$  be the computed   solution produced by the WSGD \eqref{scheme-1}-\eqref{scheme-2} or the FCD scheme \eqref{scheme-c1}-\eqref{scheme-c2} for \eqref{vp-eq-1}-\eqref{vp-b-1}  with the right hand side data $f=\mathcal{L}u,~ \mathcal{L}u^r, ~\mathcal{L}u^s$,   respectively.

Obviously, we have
\begin{equation}\label{eq-b1}
u_h(x_j)= u^r_h(x_j)+\xi^s u^s_h(x_j),\;  1\leq j\leq M-1.
\end{equation}
Moreover, derived from \eqref{eq:singularrate} and \eqref{eq:smoothrate}, for $1\leq j\leq M-1, $  we get
\begin{subequations}\label{eq-approxi-2}
\begin{eqnarray}
&&u^r(x_j)=u^r_h(x_j)+\mathcal{O}(h^2),  \label{eq-b2}\\
&&u^s(x_j)=u^s_h(x_j)+\mathcal{O}(h^{\sigma}).\label{eq-b3}
\end{eqnarray}
\end{subequations}
Thus, following from  \eqref{eq-b1}-\eqref{eq-approxi-2},  we arrive at
\begin{eqnarray}
u(x_j)&=&u^r(x_j)+\xi^su^s(x_j) \nonumber\\
& =&u^r_h(x_j)+\mathcal{O}(h^2)+\xi^s[u^s_h(x_j)+\mathcal{O}(h^{\sigma})] \nonumber\\
&=&u^r_h(x_j)+\xi^su^s_h(x_j)+\xi^s\cdot \mathcal{O}(h^{\sigma})+\mathcal{O}(h^2) \nonumber\\
&=&u_h(x_j)+\xi^s\cdot \mathcal{O}(h^{\sigma})+\mathcal{O}(h^2). \label{eq-b4}
\end{eqnarray}

From above derivation, it is readily to know that the WSGD scheme \eqref{scheme-1}-\eqref{scheme-2} and the FCD scheme \eqref{scheme-c1}-\eqref{scheme-c2} can only be of  $\sigma$-th order accuracy instead of the second-order accuracy when the solution contains a singular term, which does usually happen in physical problems albeit with smooth data.

Next, for \eqref{vp-eq-1} with the solution in the form of \eqref{split-part}, we use the extrapolation technique to improve the accuracy of schemes \eqref{scheme-1}-\eqref{scheme-2} and \eqref{scheme-c1}-\eqref{scheme-c2}. Note that
\begin{eqnarray}
&&u(x_j)=u_h(x_j)+\xi^s[u^s(x_j)-u^s_h(x_j)]+\mathcal{O}(h^2). \label{eq-b5}
\end{eqnarray}
So the corrected solution $u_h(x_j)+\xi^s[u^s(x_j)-u^s_h(x_j)]$ will dramatically improve the accuracy and convergence rate compared to the previous approximation $u_h(x_j).$
Although   it  is hard to  solve $\xi^s$ analytically,  we can compute an approximate value of $\xi^s$  with the help of extrapolation technique.

Denote $u_{h/2},~ u^r_{h/2},~ u^s_{h/2},$  as the computed   solution produced by the WSGD \eqref{scheme-1}-\eqref{scheme-2} or the FCD scheme \eqref{scheme-c1}-\eqref{scheme-c2} for \eqref{vp-eq-1}-\eqref{vp-b-1}  with the right hand side data $f=\mathcal{L}u,~ \mathcal{L}u^r, ~\mathcal{L}u^s$,   respectively. Then by \eqref{eq:singularrate}, \eqref{eq:smoothrate} and \eqref{eq-b4}, we obtain
\begin{subequations}\label{eq-approxi-1}
\begin{eqnarray}
u^r(x_j)&=&u^r_{h/2}(x_j)+\mathcal{O}(\frac{h^2}{4}),  \label{eq-b6}\\
u^s(x_j)&=&u^s_{h/2}(x_j)+\mathcal{O}(\frac{h^{\sigma}}{2^{\sigma}}), \label{eq-b7}\\
u(x_j)&=&u_{h/2}(x_j)+\mathcal{O}(\frac{h^{\sigma}}{2^{\sigma}}), \label{eq-b8}
\end{eqnarray}
\end{subequations}
for $ 1\leq j \leq M-1.$

Combining \eqref{eq-b1}-\eqref{eq-approxi-1}, we have
\begin{eqnarray}
u_{h/2}(x_j)- u_{h}(x_j)&=& u^r_{h/2}(x_j)- u^r_{h}(x_j)+\xi^s[ u^s_{h/2}(x_j)- u^s_{h}(x_{j})].
\end{eqnarray}
So it follows that
\begin{eqnarray}\label{eq-b9}
&&\xi^s=\frac{u_{h/2}(x_j)- u_{h}(x_{j})}{u^s_{h/2}(x_j)- u^s_{h}(x_{j})}-\frac{u^r_{h/2}(x_j)- u^r_{h}(x_{j})}{u^s_{h/2}(x_j)- u^s_{h}(x_{j})} .
\end{eqnarray}
Let $$\xi^s_h(x_j): =\frac{u_{h/2}(x_j)- u_{h}(x_{j})}{u^s_{h/2}(x_j)- u^s_{h}(x_{j})}.$$
Again, using \eqref{eq-b5}-\eqref{eq-approxi-1}, \eqref{eq-b9}, we get
\begin{eqnarray}
u(x_j)&=&u_h(x_j)+\bigg[\xi^s_h(x_j)-\frac{u^r_{h/2}(x_j)- u^r_{h}(x_{j})}{u^s_{h/2}(x_j)- u^s_{h}(x_{j})}\bigg][u^s(x_j)-u^s_h(x_j)]+\mathcal{O}(h^2)\nonumber\\
&=&u_h(x_j)+\xi^s_h(x_j)[u^s(x_j)-u^s_h(x_j)]-\frac{u^r_{h/2}(x_j)- u^r_{h}(x_{j})}{u^s_{h/2}(x_j)- u^s_{h}(x_{j})}[u^s(x_j)-u^s_h(x_j)]+\mathcal{O}(h^2)\nonumber\\
&=&u_h(x_j)+\xi^s_h(x_j)[u^s(x_j)-u^s_h(x_j)]+\mathcal{O}(h^2).
\end{eqnarray}

Thus, the corrected numerical solution $u(x_j)\approx u_h(x_j)+\xi^s_h(x_j)[u^s(x_j)-u^s_h(x_j)]$ is of second-order convergence. We summarize  the improved algorithm as follows.
\begin{algo}\label{algo1}
\quad

\textbf{Step 1.} Using the  WSGD scheme \eqref{scheme-1}-\eqref{scheme-2} to compute the original  problem \eqref{vp-eq-1}-\eqref{vp-b-1}, or the FCD scheme \eqref{scheme-c1}-\eqref{scheme-c2} for
solving \eqref{vp-R}, with the stepsize $h$ and  $h/2$. The numerical solutions are denoted as $u_h$ (coarse grid) and $u_{h/2}$ (fine grid) respectively.
	
	\textbf{Step 2.} Calculate the right hand data $f_s$ analytically  with the corrected  function $u_s$ given by \eqref{eq:singularterm}, which can be pre-calculated based on the fractional derivatives of power functions; see \cite{Samko1993}.  Then repeating the step 1, solve \eqref{vp-eq-1} with the right hand side function $f_s$ with the stepsize $h$ and $h/2$, and denote the numerical solutions as $u^s_{h}$ and $u^s_{h/2}$ respectively.
	
	\textbf{Step 3.}  Calculate the  strength $\xi_{h}^s$ by the  following formula
	\begin{eqnarray}
	&&\xi_{h}^s(x_j)=\frac{u_{h/2}(x_j)- u_{h}(x_{j})}{u^s_{h/2}(x_j)- u^s_{h}(x_{j})},\quad  1\leq j\leq M-1.
	\end{eqnarray}
	
	\textbf{Step 4.} Calculate the corrected solution $u^{c}_{h}=u_{h}+\xi_{h}^s(u^s-u^s_{h})$ for $  1\leq j\leq M-1.$
\end{algo}

\begin{rem}
	Readers may find that it is somewhat wasteful not to correct the solution on the fine grid. Actually, it is readily to see that the numerical solution  $u_{h/2}$ do recover  the second-order accuracy at the grid points $x_j$ with $x_j=jh$.   As to the  grid points $x_{j+\frac{1}{2}}=(x_j+x_{j+1})/2$ for $0\leq j\leq M-1,$ similarly, we have
	\begin{eqnarray*}
	&&u(x_{j+\frac{1}{2}})\nonumber\\
	 &=&u_{h/2}(x_{j+\frac{1}{2}})+\bigg[\xi^s_h(x_{j+1})-\frac{u^r_{h/2}(x_j)- u^r_{h}(x_{j})}{u^s_{h/2}(x_j)- u^s_{h}(x_{j})}\bigg][u^s(x_{j+\frac{1}{2}})-u^s_{h/2}(x_{j+\frac{1}{2}})]+\mathcal{O}(h^2)\nonumber\\
	 &=&u_{h/2}(x_{j+\frac{1}{2}})+\xi^s_{h}(x_{j+1})[u^s(x_{j+\frac{1}{2}})-u^s_{h/2}(x_{j+\frac{1}{2}})]
-\frac{u^r_{h/2}(x_j)- u^r_{h}(x_{j})}{u^s_{h/2}(x_j)- u^s_{h}(x_{j})}[u^s(x_{j+\frac{1}{2}})-u^s_{h/2}(x_{j+\frac{1}{2}})]\nonumber\\
&&+\mathcal{O}(h^2),\quad 0\leq j\leq M-2.
	\end{eqnarray*}
To apply the improved algorithm to the time-dependent case, we need the corrected value in the fine grid. For this end,  we take
    \begin{eqnarray*}
u^{c}_{h}(x_{j+\frac{1}{2}})&=&u_{h/2}(x_{j+\frac{1}{2}})+\xi^s_{h}(x_{j+1})[u^s(x_{j+\frac{1}{2}})-u^s_{h/2}(x_{j+\frac{1}{2}})],\quad 0\leq j\leq M-2,\\
u^{c}_{h}(x_{M-\frac{1}{2}})&=&u_{h/2}(x_{M-\frac{1}{2}})+\xi^s_{h}(x_{M-1})[u^s(x_{M-\frac{1}{2}})-u^s_{h/2}(x_{M-\frac{1}{2}})].
	\end{eqnarray*}
We will provide an example to show the efficiency and accuracy of the above scheme; see Example \ref{exm3}.
\end{rem}

\begin{rem}\label{rem43}
	In real  application,   the singularity of the  solution may be hierarchical and the correction of leading singular term may not ensure the convergence rate up to second-order. While the accuracy of numerical solutions can still be improved significantly by using the proposed algorithm \ref{algo1}; see Example \ref{exm1}-\ref{exm2}.  One also can reuse Algorithm \ref{algo1} to correct multi-term singularities of solution, if it has, to further enhance the accuracy and convergence rate.
\end{rem}

\section{Numerical examples}
In this section we present some numerical examples to verify our theoretical findings.

For convenience, we abbreviate the improved WSGD scheme as \textbf{I-WSGD}, the improved FCD scheme as \textbf{I-FCD}, which are obtained by applying the proposed Algorithm \ref{algo1} to the original WSGD scheme\eqref{scheme-1}-\eqref{scheme-2} and the FCD scheme  \eqref{scheme-c1}-\eqref{scheme-c2}, respectively.

In Examples \ref{exm1} and \ref{exm2}, we show the accuracy and convergence rate of the I-WSGD and I-FCD schemes for the boundary value problem \eqref{vp-eq-1}-\eqref{vp-b-1}.
We denote by  $u_j$ an approximation to $u(x_j)$ for \eqref{vp-eq-1}-\eqref{vp-b-1} obtained by the numerical schemes in the present work with spacial step size $h,$  and we measure the errors in the following sense:
\[E_\infty(h)=\max_{0\leq j\leq  M}|u^{\rm ref}_j-u_j|.
 \]
If the exact solution $u(x)$ is available, then we take $u^{\rm ref}_j=u(x_j)$; otherwise we compute the reference solution  $u^{\rm ref}_j$ with the step size $h=2^{-15}$. In Example \ref{exm3}, the I-WSGD scheme  will be applied to a time-dependent SpFDE.

	\begin{example}\label{exm1}
		Consider the problem \eqref{vp-eq-1}-\eqref{vp-b-1} with the left-sided fractional derivative, that is $\theta=1.$  We take $a=0$, $b=1$ and $\alpha=1.$
		
		{Case I} \quad Choose suitable $f$ such that the exact solution to \eqref{vp-eq-1}-\eqref{vp-b-1} is  $u(x)=(x^2+x^{\beta +1}+x^{\beta -1})(1-x).$
		
		{Case II} \quad Take the corresponding right hand side function in \eqref{vp-eq-1}-\eqref{vp-b-1} as $f(x)=x+1$.
	\end{example}
	For both Case I and Case II, we take the leading weak singular term as $u^s=x^{\beta -1}(1-x)$, then the corresponding right hand side function is  $f^s=x^{\beta -1}(1-x)+\Gamma(\beta +1).$  Tables \ref{table1} and \ref{table2} show that both for the case with a given non-smooth solution (Case I), and for the case where the exact solution is unknown (Case II),  the accuracy and convergence rate of numerical solutions from the I-WSGD scheme are improved significantly, compared to the numerical solution from  the original WSGD scheme \eqref{scheme-1}-\eqref{scheme-2}, of which the convergence order is only $\beta-1$. It is worth to  mention that, for $\beta=1.1$, data in Table \ref{table1} illustrate that the I-WSGD scheme is of second-order convergence, while the convergence rate is only round one shown in Table \ref{table2}. It implies that to get second-order accuracy for small $\beta$ (close to one), we need to consider more singularities in the improved scheme than the leading weak singularity; see also Remark  \ref{rem43}.

\begin{figure}[!h]
	\centering
	 \includegraphics[width=0.33\linewidth]{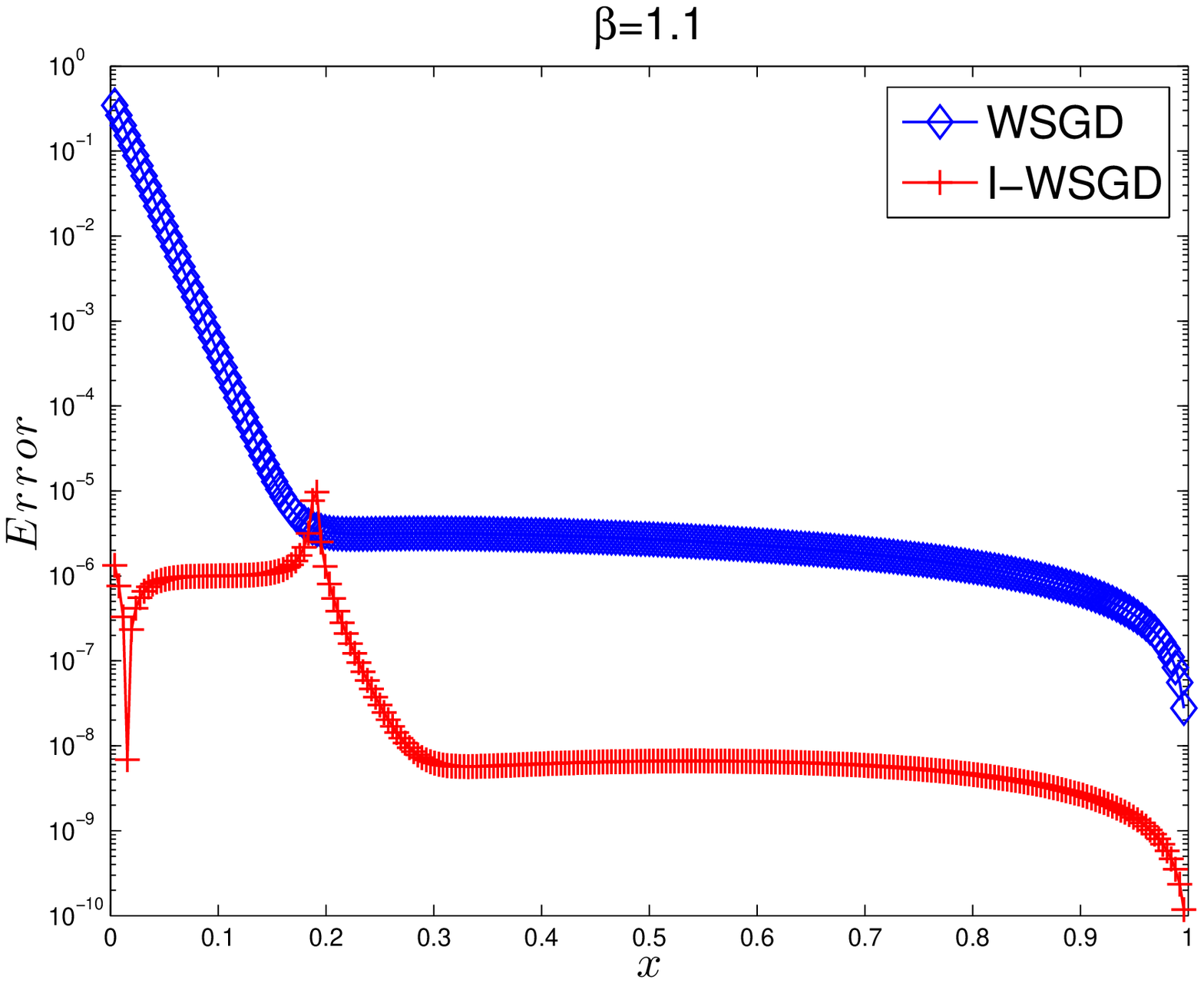}~\includegraphics[width=0.33\linewidth]{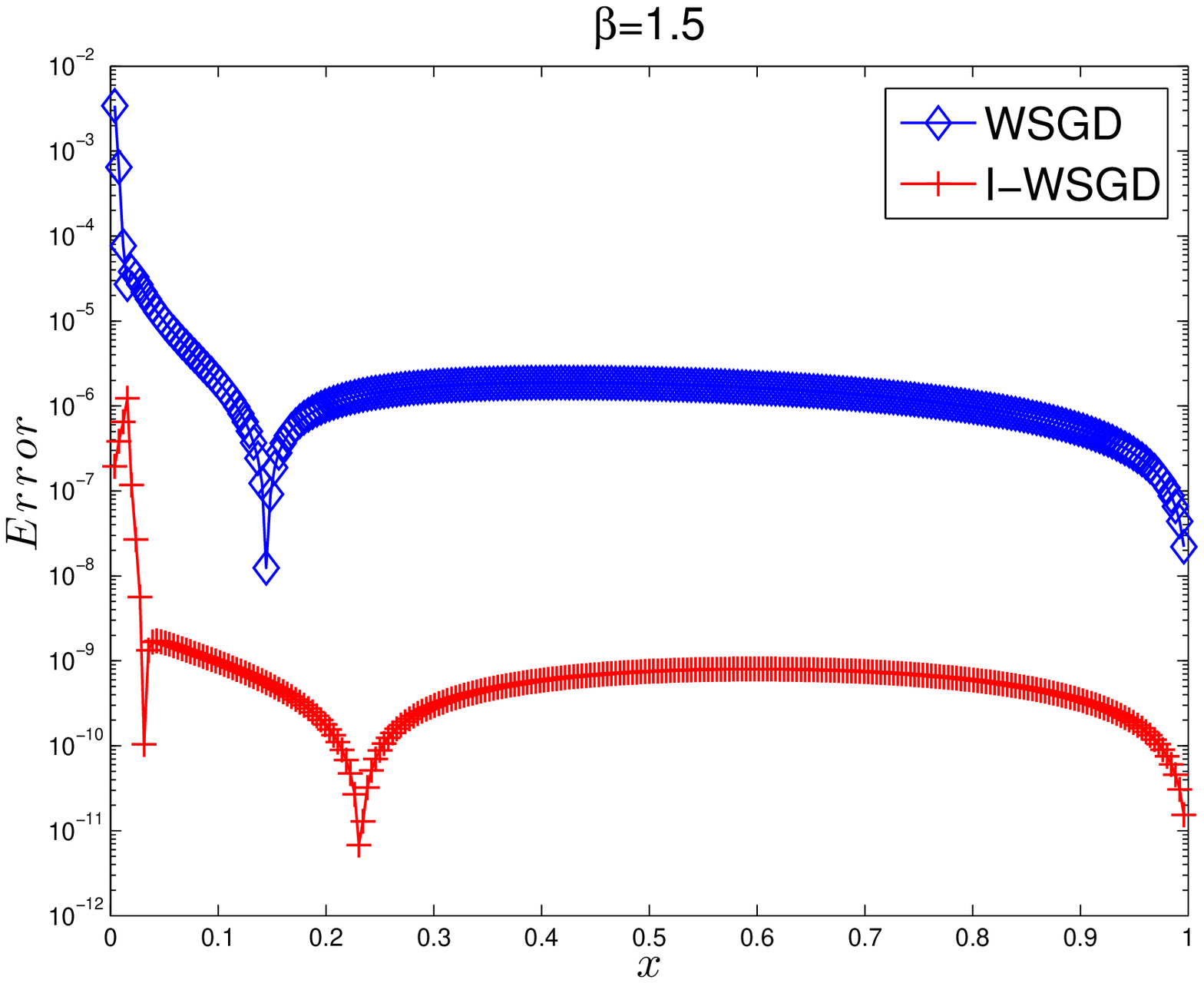}~\includegraphics[width=0.33\linewidth]{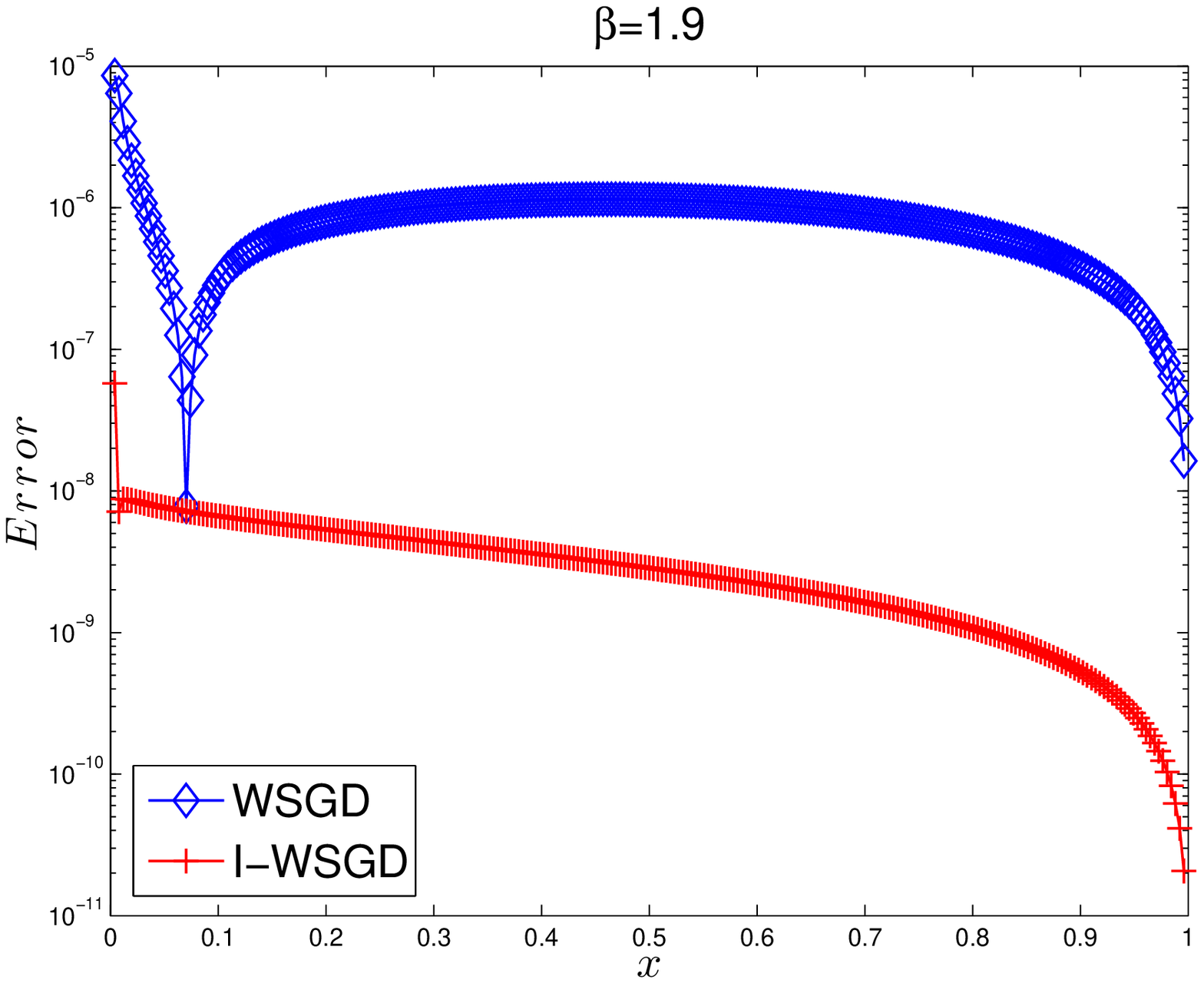}
	\caption{Comparison of point-wise errors for the WSGD   scheme \eqref{scheme-1}-\eqref{scheme-2} and the I-WSGD scheme  for \eqref{vp-eq-1}-\eqref{vp-b-1} with left sided fractional derivative (Example \ref{exm1}, Case I). $h=2^{-9}$.}
	\label{fig:solleft}
\end{figure}

\begin{figure}[!h]
	\centering
	 \includegraphics[width=0.33\linewidth]{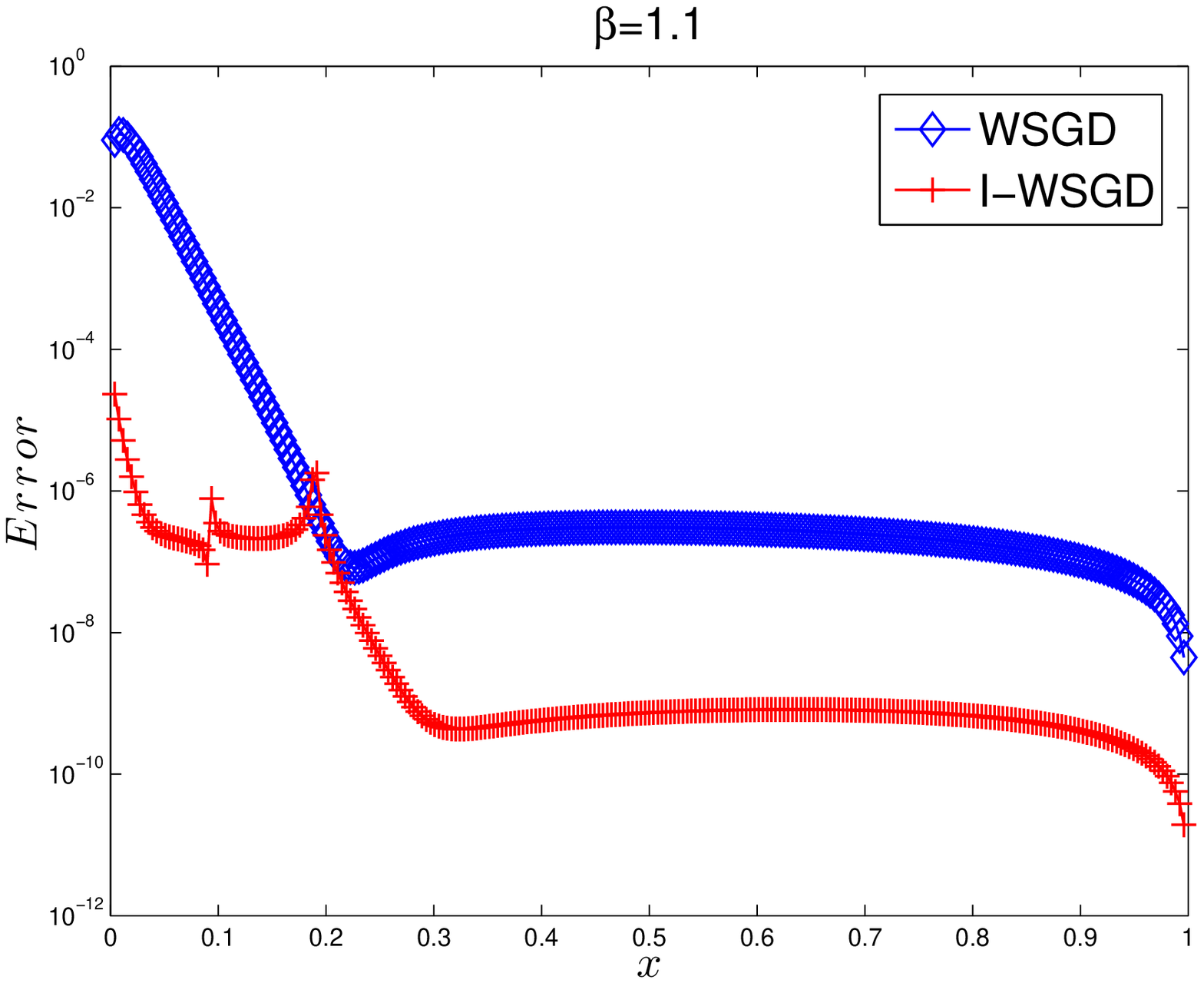}~\includegraphics[width=0.33\linewidth]{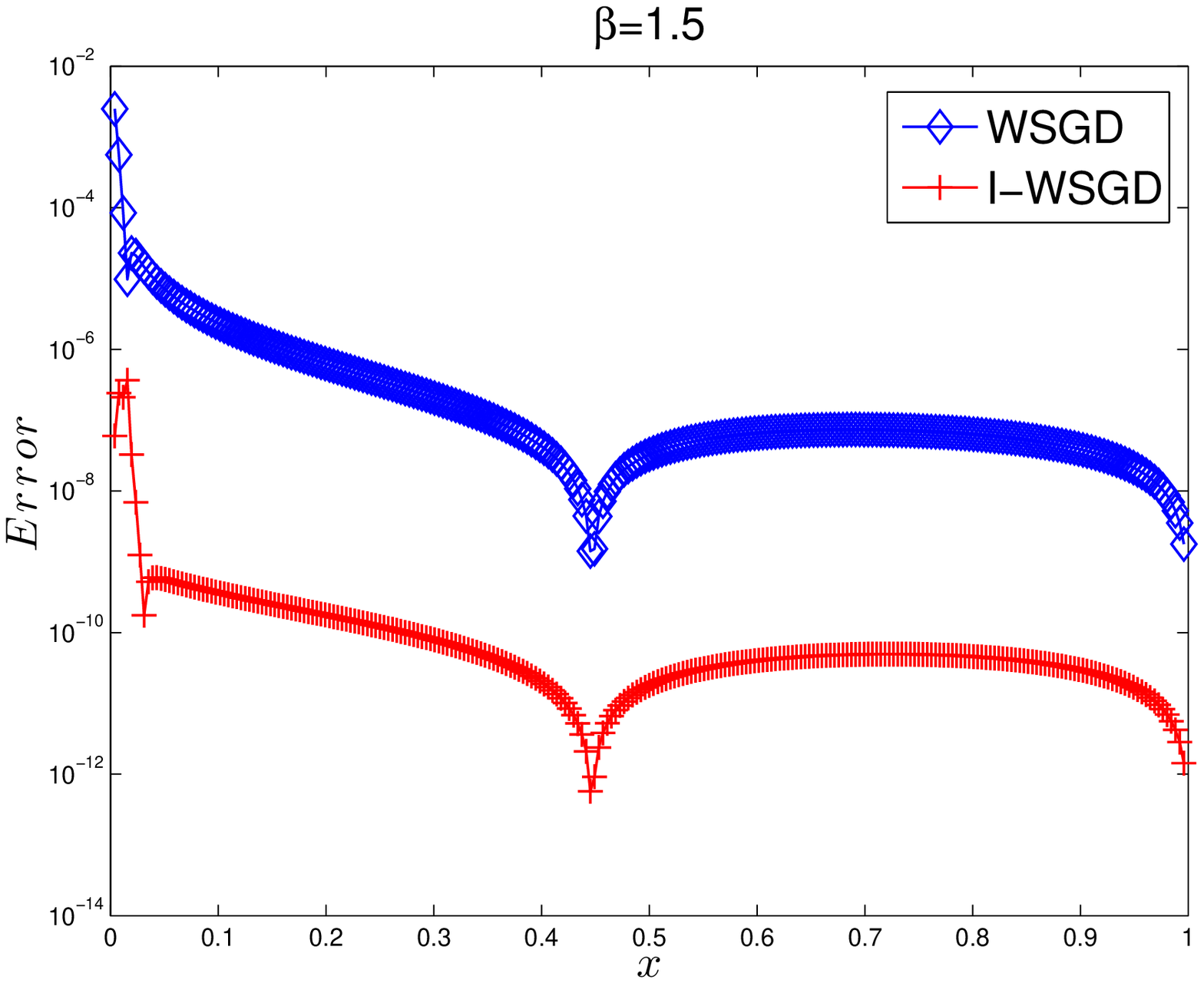}~\includegraphics[width=0.33\linewidth]{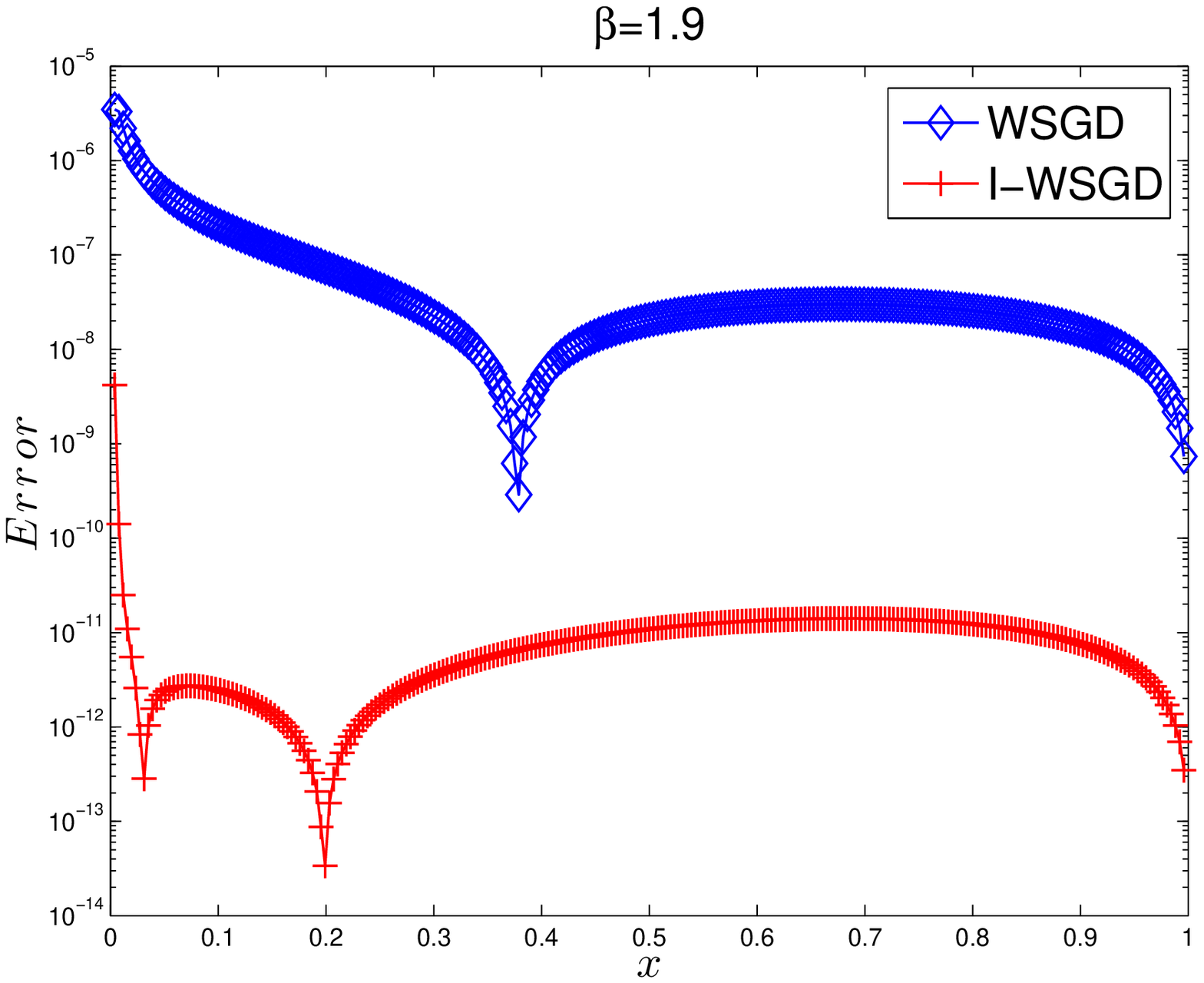}
	\caption{Comparison of point-wise errors for the WSGD   scheme \eqref{scheme-1}-\eqref{scheme-2} and the I-WSGD scheme  for \eqref{vp-eq-1}-\eqref{vp-b-1} with left-sided fractional derivative (Example \ref{exm1}, Case II). $h=2^{-10}$. }
	\label{fig:nosolleft}
\end{figure}

	In Figs \ref{fig:solleft} and \ref{fig:nosolleft}, we show the behavior of point-wise errors of numerical solutions from the I-WSGD scheme and the original WSGD scheme for the problem \eqref{vp-eq-1}-\eqref{vp-b-1} with non-smooth solution. Figs \ref{fig:solleft} and \ref{fig:nosolleft} illustrates that the error has been reduced remarkably by applying the Algorithm \ref{algo1}, even for Case II where we do not know what the exact solution is.

\begin{table}[!htb]
	\centering
	\caption{\scriptsize Comparison of accuracy and convergence rate between the original WSGD scheme \eqref{scheme-1}-\eqref{scheme-2} and the I-WSGD scheme for \eqref{vp-eq-1}-\eqref{vp-b-1} with left-sided fractional derivative (Example \ref{exm1}, Case I). The exact solution  is chosen as $(x^2+x^{\beta +1}+x^{\beta -1})(1-x)$.}\label{table1}
	{\scriptsize\begin{tabular}{c|cccc|cccc}
			\hline
			\multirow{2}{*}{$\beta $}& \multicolumn{4}{c|}{the WSGD scheme}  &\multicolumn{4}{c}{the  I-WSGD  scheme } \\
			\cline{2-5}
			\cline{6-9}
			& $M$        & $E_\infty(h)$   & $Rate$     & CPU time (s) & $M$     & $E_\infty(h)$   & $Rate$     & CPU time (s)  \\
			\hline
$\beta=1.1$&	512        & 4.03e-01     &          &  0.88  & 64        & 2.45e-04     &         &  0.03      \\
		   &	1024       & 3.77e-01     &     0.10  &  0.17  & 128       & 1.16e-04     &     1.08 &  0.06       \\
		   &	2048       & 3.52e-01     &     0.10  &  0.60  & 256       & 5.30e-05     &     1.13 &  0.10       \\
		   &	4096       & 3.28e-01     &     0.10  &  3.14  & 512       & 9.78e-06     &     2.44 &  0.32        \\
			\hline
$\beta=1.5$&	512        & 9.52e-03     &          &  0.05  & 64        & 1.32e-04     &         &  0.02      \\
		   &	1024       & 6.73e-03     &     0.50  &  0.15  & 128       & 2.82e-05     &     2.23 &  0.07       \\
		   &	2048       & 4.76e-03     &     0.50  &  0.52  & 256       & 6.27e-06     &     2.17 &  0.12       \\
		   &	4096       & 3.37e-03     &     0.50  &  2.34  & 512       & 1.42e-06     &     2.14 &  0.36        \\
			\hline
$\beta=1.9$&	512        & 7.44e-05     &          &  0.06  & 64        & 1.19e-05     &         &  0.02      \\
		   &	1024       & 3.99e-05     &     0.90  &  0.17  & 128       & 2.49e-06     &     2.25 &  0.06       \\
		   &	2048       & 2.14e-05     &     0.90  &  0.64  & 256       & 5.50e-07     &     2.18 &  0.12       \\
		   &	4096       & 1.15e-05     &     0.90  &  2.51  & 512       & 1.27e-07     &     2.11 &  0.33        \\
			\hline
		\end{tabular}}
	\end{table}
	
\begin{table}[!htb]
	\centering
	\caption{\scriptsize Comparison of accuracy and convergence rate between the original WSGD scheme \eqref{scheme-1}-\eqref{scheme-2} and the I-WSGD scheme  for solving \eqref{vp-eq-1}-\eqref{vp-b-1} with left-sided fractional derivative (Example \ref{exm1}, Case II). The right hand side function  is chosen as $f(x)=x+1$.}\label{table2}
	{\scriptsize\begin{tabular}{c|cccc|cccc}
			\hline
			\multirow{2}{*}{$\beta $}& \multicolumn{4}{c|}{the WSGD scheme }  &\multicolumn{4}{c}{ the I-WSGD scheme } \\
			\cline{2-5}
			\cline{6-9}
			& $M$        & $E_\infty(h)$   & $Rate$     & CPU time (s) & $M$     & $E_\infty(h)$   & $Rate$     & CPU time (s)  \\
			\hline
$\beta=1.1$&	1024       & 6.60e-01     &         &  0.08  & 64        & 5.74e-04     &        &  0.01      \\
&	2048       & 6.17e-01     &     0.10 &  0.34  & 128       & 2.68e-04     &     1.10 &  0.02       \\
&	4096       & 5.75e-01     &     0.10 &  2.05  & 256       & 1.20e-04     &     1.15 &  0.05       \\
&	8192       & 5.37e-01     &     0.10 &  15.7  & 512       & 5.31e-05     &     1.18 &  0.19        \\
			\hline
$\beta=1.5$			&	1024       & 1.01e-02     &           &  0.08  & 64        & 8.84e-05     &         &  0.01      \\
			&	2048       & 7.11e-03     &     0.50   &  0.37  & 128       & 1.86e-05     &     2.25 &  0.02       \\
			&	4096       & 5.02e-03     &     0.50   &  2.17  & 256       & 4.83e-06     &     1.95 &  0.04       \\
			&	8192       & 3.55e-03     &     0.50   &  16.6  & 512       & 1.32e-06     &     1.87 &  0.19        \\
			\hline
$\beta=1.9$			&	1024       & 4.63e-05     &          &  0.08  & 64        & 1.84e-06     &         &  0.01      \\
			&	2048       & 2.48e-05     &     0.90  &  0.39  & 128       & 3.69e-07     &     2.32 &  0.02       \\
			&	4096       & 1.33e-05     &     0.90  &  2.37  & 256       & 7.79e-08     &     2.24 &  0.06       \\
			&	8192       & 7.11e-06     &     0.90  &  16.8  & 512       & 1.75e-08     &     2.15 &  0.18        \\
\hline
		\end{tabular}}
	\end{table}


\begin{example}\label{exm2}
Consider the problem \eqref{vp-eq-1}-\eqref{vp-b-1} with the symmetrical two-sided fractional derivatives, that is $\theta=1/2$. We take $a=0$, $b=1$ and $\alpha=1.$

{Case I} \quad Choose suitable $f$ such that the exact solution to \eqref{vp-eq-1}-\eqref{vp-b-1} is  $u(x)=x^2(1-x)^2+2x^{\beta /2}(1-x)^{\beta /2}.$
	
{Case II} \quad Take the corresponding right hand side function in \eqref{vp-eq-1}-\eqref{vp-b-1} as $f(x)=1$.
\end{example}
In this example, the weak singular function is chosen as $u^s=x^{\beta /2}(1-x)^{\beta /2}$ for both Case I and Case II,  and the corresponding right side function is $f^s=x^{\beta /2}(1-x)^{\beta /2}-\cos(\beta /2 \pi)\Gamma(\beta +1).$
Similar to Example \ref{exm1}, Tables \ref{table3}-\ref{table4} show that for the equation \eqref{vp-eq-1}-\eqref{vp-b-1} with symmetric two-sided fractional derivative, the I-FCD scheme can enhance the accuracy and convergence order of numerical solutions greatly. Moreover, compared with the low convergence rate $\beta/2$ of the FCD scheme, the convergence rate of the I-FCD scheme is more than $1.5$ for $\beta=1.1$, and second-order accuracy can be obtained for $\beta=1.5$ and $1.9$.  We also use the scheme \eqref{scheme-1}-\eqref{scheme-2} to solve this example, the numerical results are similar, which we do not present here.

\begin{figure}[!h]
	\centering
	 \includegraphics[width=0.33\linewidth]{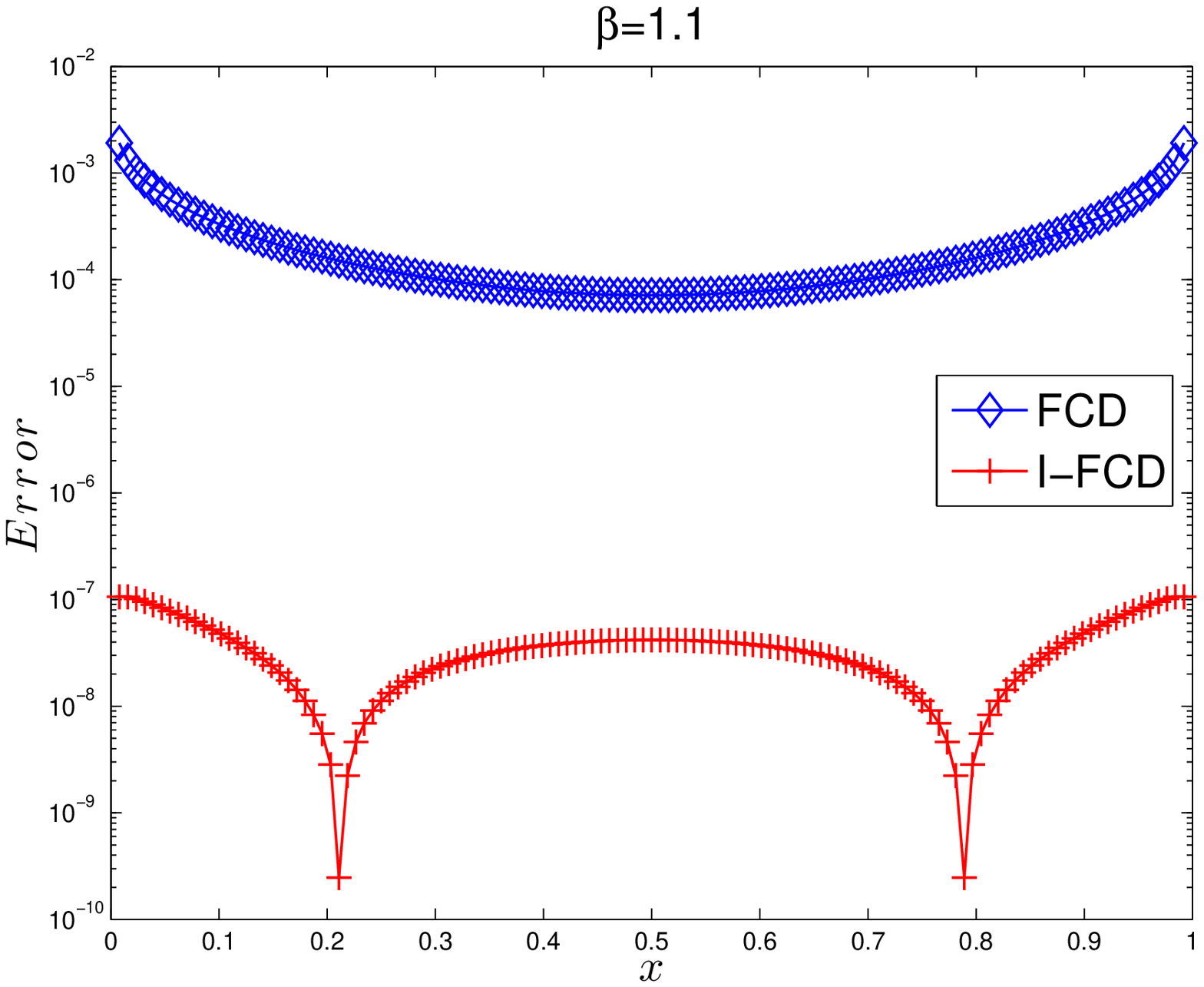}~\includegraphics[width=0.33\linewidth]{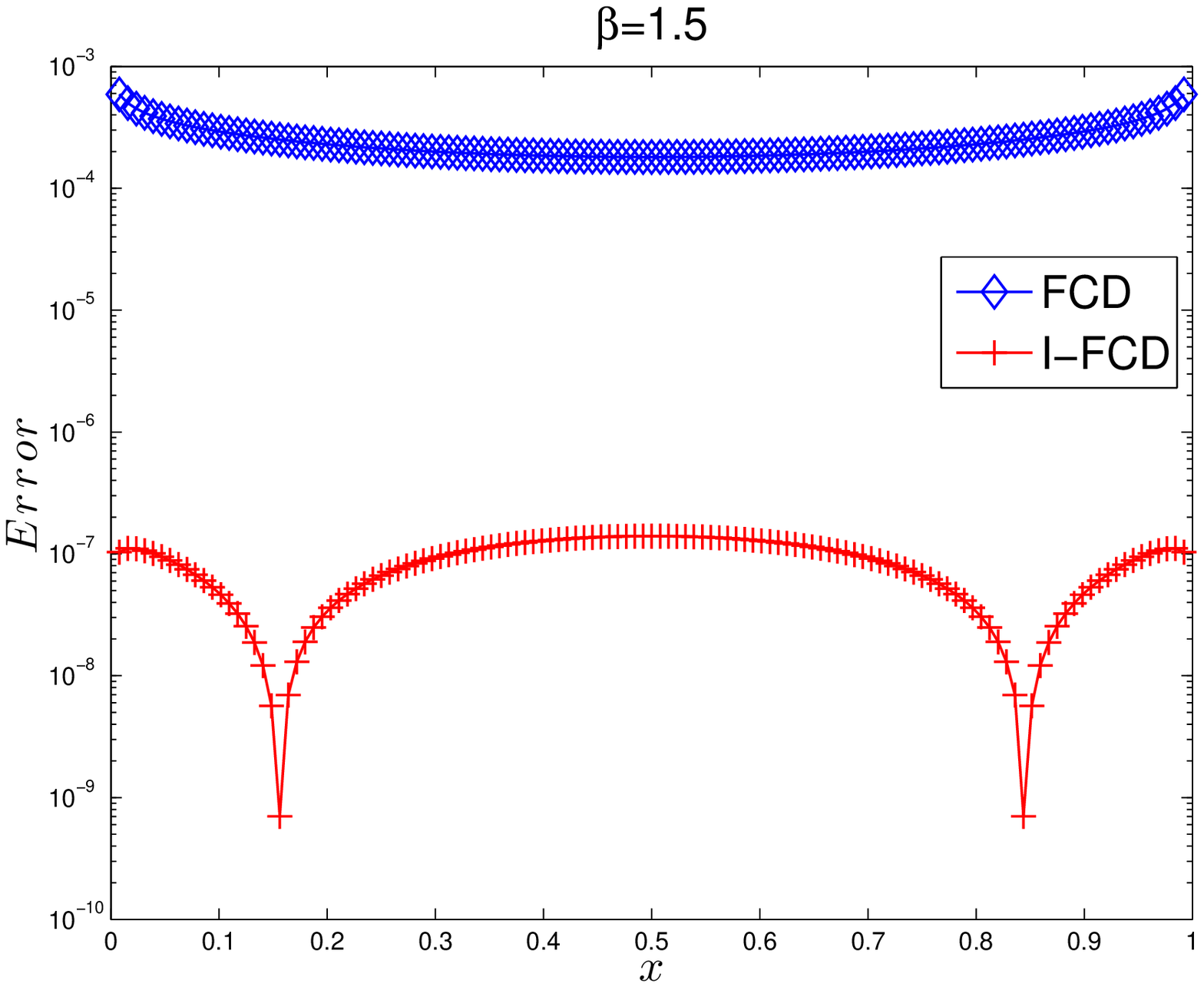}~\includegraphics[width=0.33\linewidth]{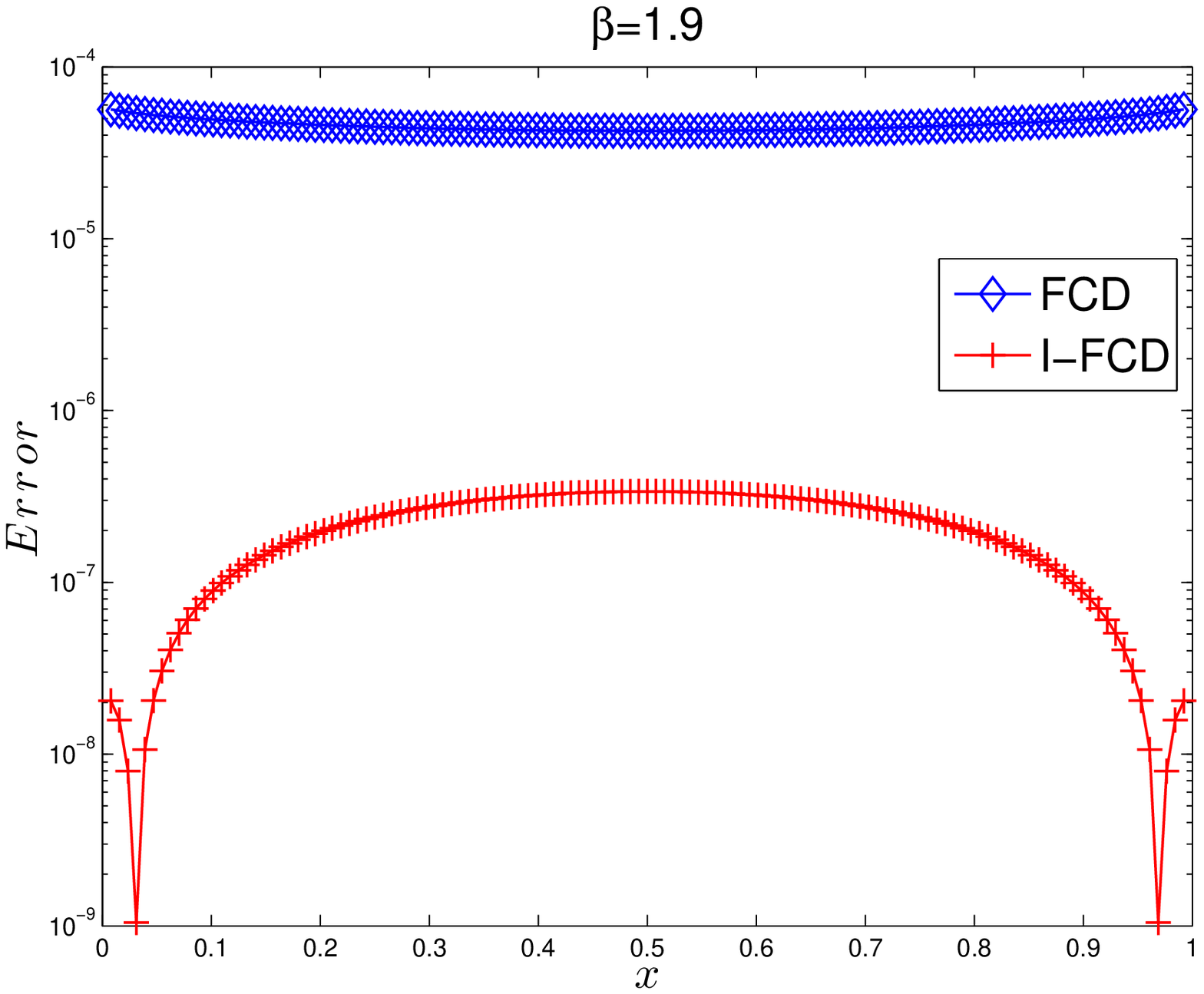}
	\caption{Comparison of point-wise errors for the   FCD scheme \eqref{scheme-c1}-\eqref{scheme-c2} and the I-FCD scheme for \eqref{vp-eq-1}-\eqref{vp-b-1} with symmetric two-sided  fractional derivative (Example \ref{exm2}, Case I). $h=2^{-9}$.}
	\label{fig:solsym}
\end{figure}

\begin{figure}[!h]
	\centering
	 \includegraphics[width=0.33\linewidth]{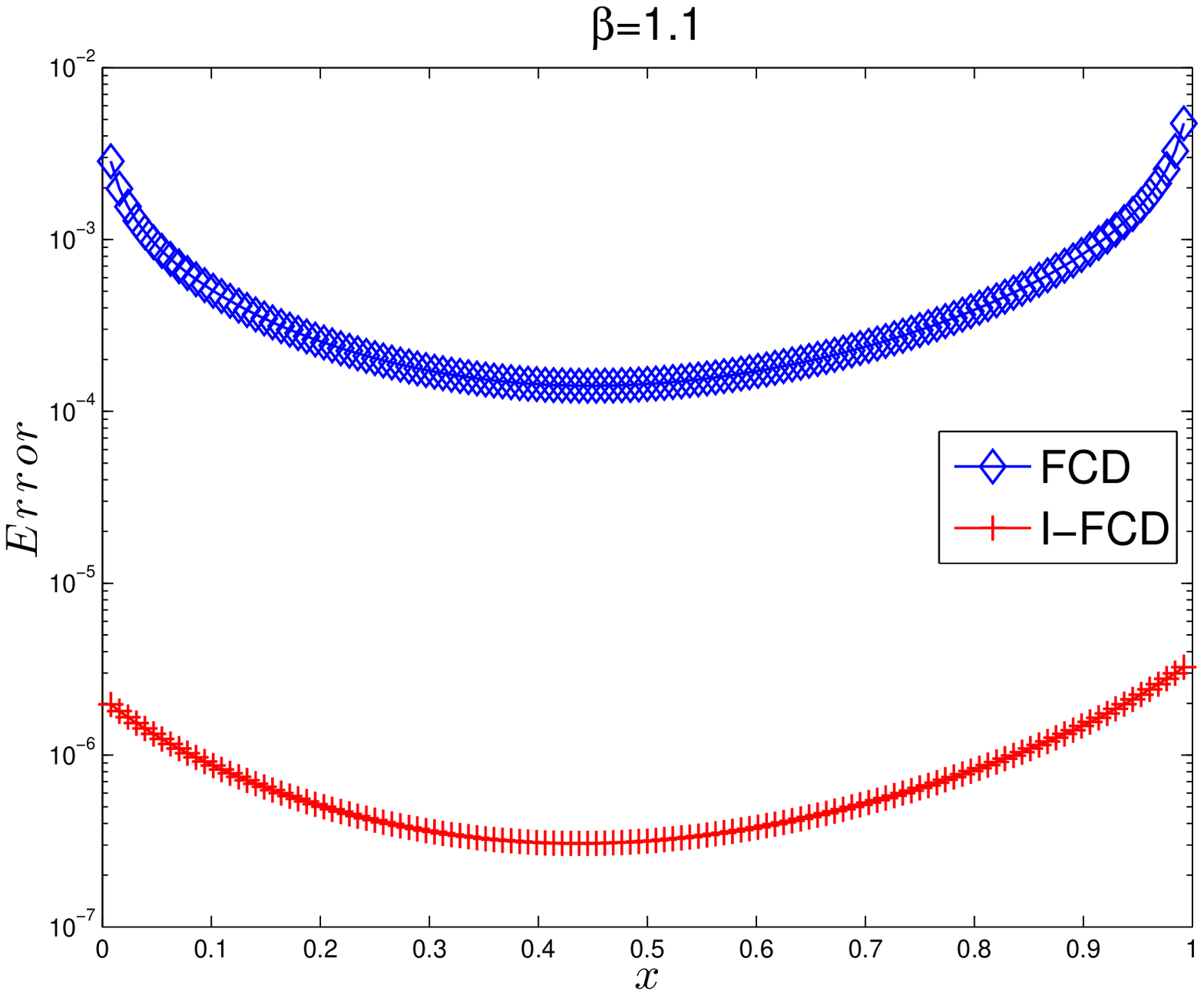}~\includegraphics[width=0.33\linewidth]{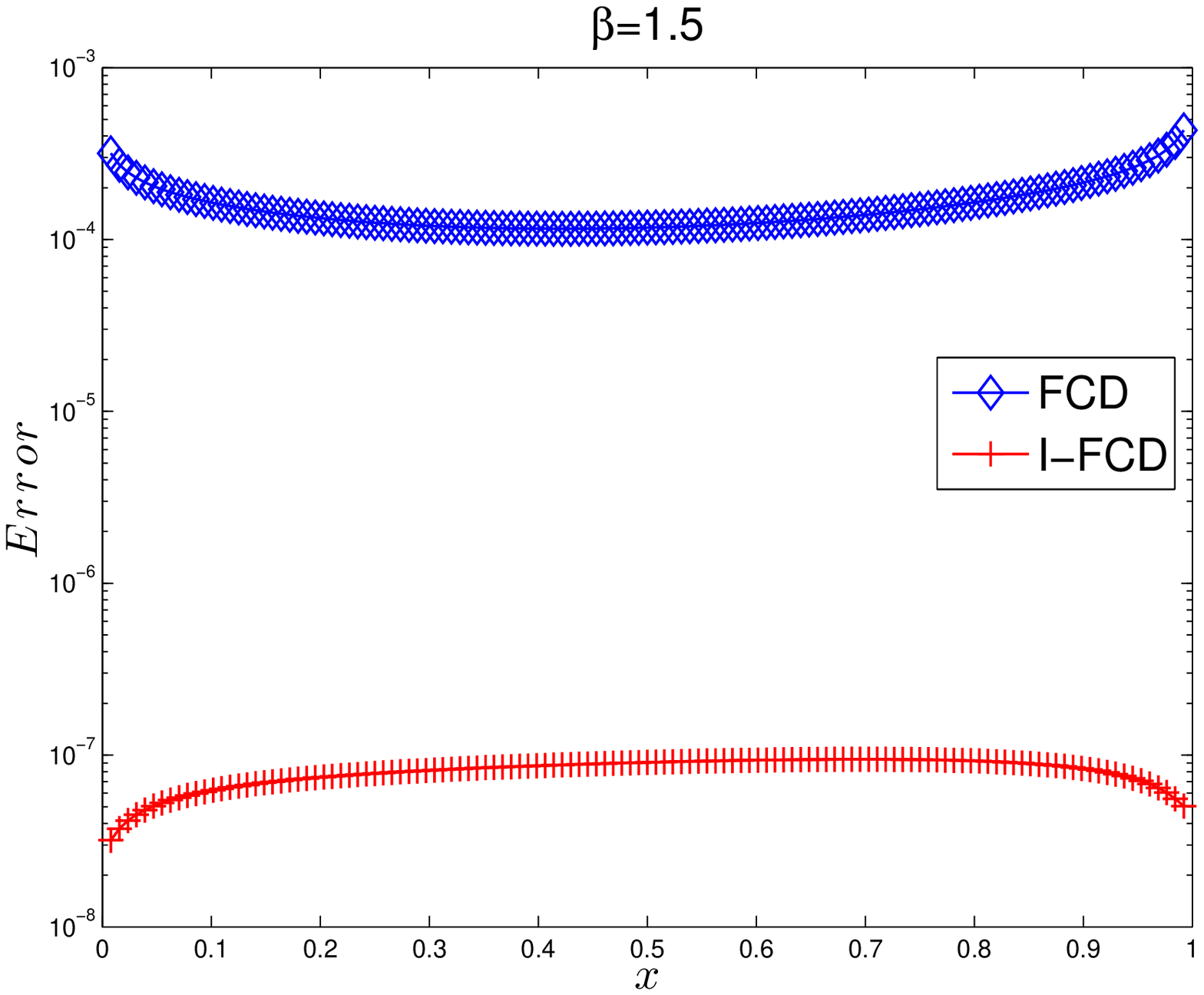}~\includegraphics[width=0.33\linewidth]{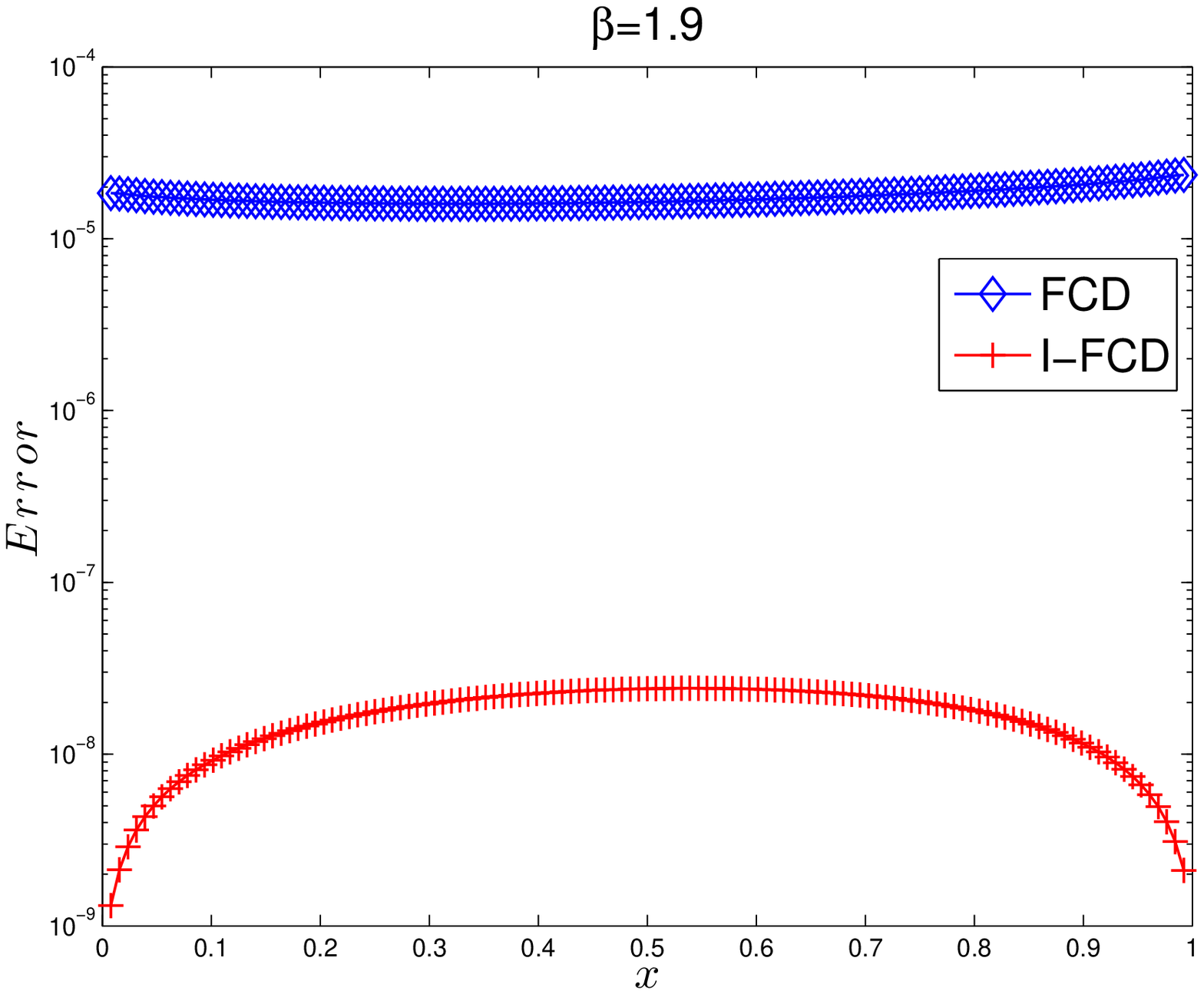}
	\caption{Comparison of point-wise errors for the FCD  scheme \eqref{scheme-c1}-\eqref{scheme-c2} and the I-FCD scheme  for \eqref{vp-eq-1}-\eqref{vp-b-1} with symmetric two-sided fractional derivative (Example \ref{exm2}, Case II). $h=2^{-10}$. }
	\label{fig:nosolsym}
\end{figure}

Further, we give Figs \ref{fig:solsym} and \ref{fig:nosolsym} to show the behavior of point-wise errors for the FCD and I-FCD schemes to solve \eqref{vp-eq-1}-\eqref{vp-b-1} with symmetric two-sided fractional derivative. It is illustrated that,  for different $\beta$, the numerical solutions from I-FCD scheme get higher accuracy than that from the original FCD scheme \eqref{scheme-c1}-\eqref{scheme-c2}.

\begin{table}[!htb]
	\centering
	\caption{\scriptsize Comparison of accuracy and convergence rate between the original FCD scheme \eqref{scheme-c1}-\eqref{scheme-c2} and the I-FCD scheme for solving \eqref{vp-eq-1}-\eqref{vp-b-1} with symmetric two-sided fractional derivative (Example \ref{exm2}, Case I). The exact solution  is chosen as $u(x)=x^2(1-x)^2+2x^{\beta /2}(1-x)^{\beta /2}$.}\label{table3}
	{\scriptsize\begin{tabular}{c|cccc|cccc}
			\hline
			\multirow{2}{*}{$\beta $}& \multicolumn{4}{c|}{the FCD scheme}  &\multicolumn{4}{c}{the I-FCD scheme} \\
			\cline{2-5}
			\cline{6-9}
			& $M$        & $E_\infty(h)$   & $Rate$     & CPU time (s) & $M$     & $E_\infty(h)$   & $Rate$     & CPU time (s)  \\
			\hline
$\beta=1.1$&	512        & 3.50e-03     &        &  0.05  & 64        & 4.18e-06     &         &  0.04      \\
&	1024       & 2.42e-03     &     0.53&  0.12  & 128       & 1.30e-06     &     1.68 &  0.06       \\
&	2048       & 1.66e-03     &     0.54&  0.42  & 256       & 3.81e-07     &     1.77 &  0.10       \\
&	4096       & 1.14e-03     &     0.54&  1.61  & 512       & 1.07e-07     &     1.83 &  0.27        \\
			\hline
$\beta=1.5$&	512        & 7.50e-04     &         &  0.04  & 64        & 1.06e-05     &         &  0.02      \\
&	1024       & 4.47e-04     &     0.75 &  0.12  & 128       & 2.49e-06     &     2.09 &  0.04       \\
&	2048       & 2.66e-04     &     0.75 &  0.38  & 256       & 5.89e-07     &     2.08 &  0.09       \\
&	4096       & 1.58e-04     &     0.75 &  1.50  & 512       & 1.40e-07     &     2.07 &  0.24        \\
\hline
$\beta=1.9$&	512        & 5.66e-05     &         &  0.04  & 64        & 2.32e-05     &          &  0.02       \\
&	1024       & 2.94e-05     &     0.95 &  0.12  & 128       & 5.65e-06     &     2.03 &  0.04        \\
&	2048       & 1.52e-05     &     0.95 &  0.36  & 256       & 1.38e-06     &     2.03 &  0.09        \\
&	4096       & 7.87e-06     &     0.95 &  1.62  & 512       & 3.38e-07     &     2.03 &  0.25         \\
\hline
		\end{tabular}}
	\end{table}

\begin{table}[!htb]
	\centering
	\caption{\scriptsize Comparison of accuracy and convergence rate  between  the original FCD scheme  \eqref{scheme-c1}-\eqref{scheme-c2} and  the I-FCD scheme for solving \eqref{vp-eq-1}-\eqref{vp-b-1} with symmetric two-sided fractional derivative (Example \ref{exm2}, Case I). The right hand side function is chosen as $f(x)=1$.}\label{table4}
	{\scriptsize\begin{tabular}{c|cccc|cccc}
			\hline
			\multirow{2}{*}{$\beta $}& \multicolumn{4}{c|}{the FCD scheme }  &\multicolumn{4}{c}{the I-FCD scheme} \\
			\cline{2-5}
			\cline{6-9}
			& $M$        & $E_\infty(h)$   & $Rate$     & CPU time (s) & $M$     & $E_\infty(h)$   & $Rate$     & CPU time (s)  \\
			\hline
$\beta=1.1$			& 1024       & 8.22e-03     &         &  0.05  & 64        & 3.37e-04     &     &  0.01      \\
			& 2048       & 5.68e-03     &     0.53 &  0.38  & 128       & 1.30e-04     &     1.37 &  0.02       \\
			& 4096       & 3.90e-03     &     0.54 &  1.58  & 256       & 4.59e-05     &     1.50 &  0.05       \\
			& 8192       & 2.68e-03     &     0.54 &  10.8  & 512       & 1.55e-05     &     1.57 &  0.17        \\
			\hline
$\beta=1.5$			&	1024       & 5.16e-04     &         &  0.05  & 64        & 2.17e-05     &         &  0.01      \\
			&	2048       & 3.07e-04     &     0.75 &  0.27  & 128       & 5.67e-06     &     1.94 &  0.02       \\
			&	4096       & 1.83e-04     &     0.75 &  1.44  & 256       & 1.46e-06     &     1.96 &  0.03       \\
			&	8192       & 1.09e-04     &     0.75 &  10.8  & 512       & 3.73e-07     &     1.97 &  0.13        \\
			\hline
$\beta=1.9$			&	1024       & 2.35e-05     &         &  0.05  & 64        & 6.22e-06     &     &  0.01      \\
			&	2048       & 1.22e-05     &     0.95 &  0.44  & 128       & 1.55e-06     &     2.00 &  0.02       \\
			&	4096       & 6.33e-06     &     0.95 &  1.42  & 256       & 3.87e-07     &     2.00 &  0.04       \\
			&	8192       & 3.28e-06     &     0.95 &  10.7  & 512       & 9.67e-08     &     2.00 &  0.14        \\
			\hline
		\end{tabular}}
	\end{table}		

				
\begin{example}\label{exm3}
In this example,  we consider the following SpFDE with the left-sided fractional derivative:
\begin{eqnarray}\label{td}
u_t(x,t)=\ _0D_x^{\beta }u(x,t)+f(x,t),  \quad 0<x<1, \quad 0<t<1,
\end{eqnarray}
which admits the solution
%$$u(x,t)=[x^{\beta /2}(1-x)^{\beta /2}+x^2(1-x)^2]t^3.  $$
$$u(x,t)=(x^{\beta-1}+x^2+x^{1+\beta})(1-x)t^3.  $$
\end{example}
For an integer $N,$ take $\tau=1/N,$ and denote $t_n=n\tau,$  and $t_{n-\frac{1}{2}}=(t_n+t_{n-1})/2$ for $n=1,2,\ldots, N$.   In terms of  temporal discretization, we adopt the classical  second-order Crank-Nikolson scheme. The semi-discretized scheme  reads
\begin{eqnarray}
u(x_j,t_n)- \frac{\tau}{2}\ _0D_x^{\beta }u(x_j,t_n) = u(x_j,t_{n-1})+\frac{\tau}{2}\ _0D_x^{\beta }u(x_j,t_{n-1})+\tau f(x_j,t_{n-\frac{1}{2}})+\mathcal{O}(\tau^3).
\end{eqnarray}
Let $u_j^n$ be the approximation of $u(x_j,t_n)$ at the grid point $ (x_j, t_n).$  We use the WSGD scheme \eqref{scheme-1}-\eqref{scheme-2} to discretize the spatial fractional derivative.  Then the fully discretized difference scheme (\textbf{CN-WSGD}) is given as follows.
\begin{eqnarray}\label{td-s}
u_j^n-\frac{\tau}{2}\delta_{x,-}^{\beta }u_j^n= u_j^{n-1}+\frac{\tau}{2}\delta_{x,-}^{\beta } u_j^{n-1}+\tau f(x_j,t_{n-\frac{1}{2}}),\quad  1\leq j\leq M-1,\quad 1\leq n\leq N.
\end{eqnarray}
To apply the Algorithm \ref{algo1}, in each time step,  we take the same leading weak singular term as that in Example \ref{exm1}. Then the numerical solution $u_j^n$ produced by the I-WSGD scheme is updated in next time step. For convenience, we denote the scheme \eqref{td-s} with the use of the Algorithm \ref{algo1} as the \textbf{CN-I-WSGD} scheme.
				
We measure errors and convergence rate in this example as follows. Let $$E_{\infty}^N(h)=\max_{1\leq j\leq M-1}|u(x_j,t_N)-u_{j}^N|$$
and assume
$$E_{\infty}^N(h) = O(h ^p)+O(\tau^q).$$ If $\tau$ is sufficiently small, then $E_{\infty}^N(h) \approx O(h^p)$. Consequently,
$\frac{E_{\infty}^N(2h)}{E_{\infty}^N(h)}\approx 2^p$ and  $p \approx  \log_2\left(\frac{E_{\infty}^N(2h)}{E_{\infty}^N(h)}\right)$ is the convergence rate with respect to the spatial step size.
				
\begin{table}[!htb]
	\centering
	\caption{\scriptsize Comparison of accuracy and convergence rate between the CN-WSGD scheme \eqref{td-s} and the improved scheme CN-I-WSGD  for solving the SpFDE \eqref{td} with left-sided fractional derivative (Example \ref{exm3}). The exact solution  is chosen as $u(x,t)=(x^{\beta-1}+x^2+x^{1+\beta})(1-x)t^3$ and $\tau=10^{-3}.$}\label{table5}
	{\scriptsize\begin{tabular}{c|cccc|cccc}
			\hline
			\multirow{2}{*}{$\beta $}& \multicolumn{4}{c|}{the CN-WSGD scheme}  &\multicolumn{4}{c}{the CN-I-WSGD scheme} \\
			\cline{2-5}
			\cline{6-9}
			& $M$        & $E^N_\infty(h)$   & $Rate$     & CPU time (s) & $M$     & $E^N_\infty(h)$   & $Rate$     & CPU time (s)  \\
			\hline
$\beta=1.4$&	16       & 9.33e-02     &          &  1.55 & 4       & 7.84e-03     &         &  1.29      \\
&	32       & 7.43e-02     &     0.33  &  2.88 & 8       & 1.31e-03     &     2.58 &  2.16       \\
&	64       & 5.73e-02     &     0.38  &  5.70 & 16      & 3.68e-04     &     1.83 &  3.99       \\
&	128      & 4.37e-02     &     0.39  &  11.2 & 32      & 8.72e-05     &     2.08 &  7.63        \\
			\hline
$\beta=1.8$&	16       & 4.81e-03     &           &  1.44  & 4       & 2.11e-03     &          &  1.30      \\
&	32       & 3.04e-03     &     0.66   &  2.75  & 8       & 3.74e-04     &     2.50  &  2.19       \\
&	64       & 1.79e-03     &     0.76   &  5.65  & 16       & 6.50e-05     &     2.52 &  3.93       \\
&	128      & 1.04e-03     &     0.79   &  11.4  & 32       & 1.53e-05     &     2.09 &  7.64        \\
			\hline
		\end{tabular}}
	\end{table}
					
\begin{figure}[!h]
	\centering
	 \includegraphics[width=0.45\linewidth]{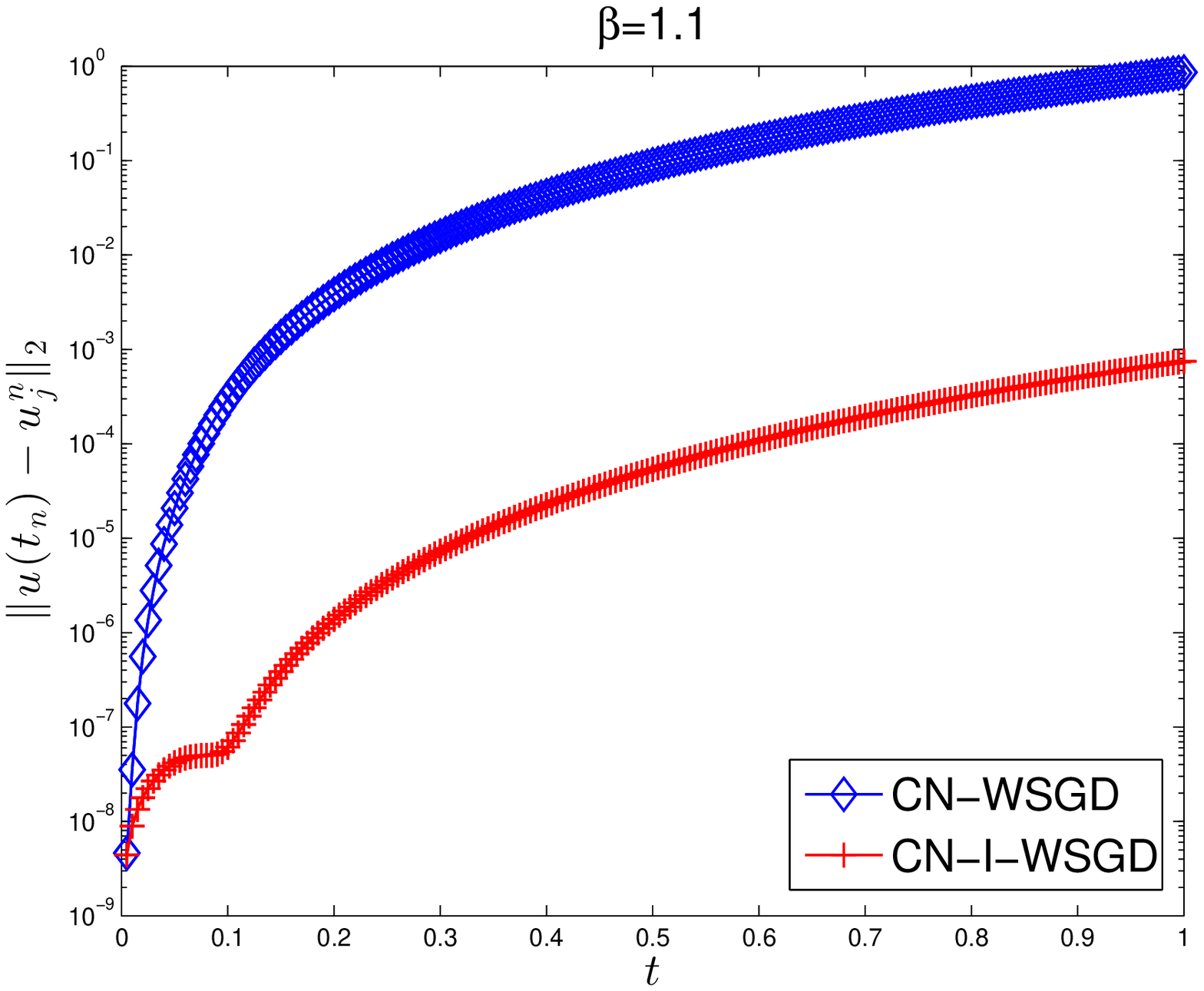}~\includegraphics[width=0.45\linewidth]{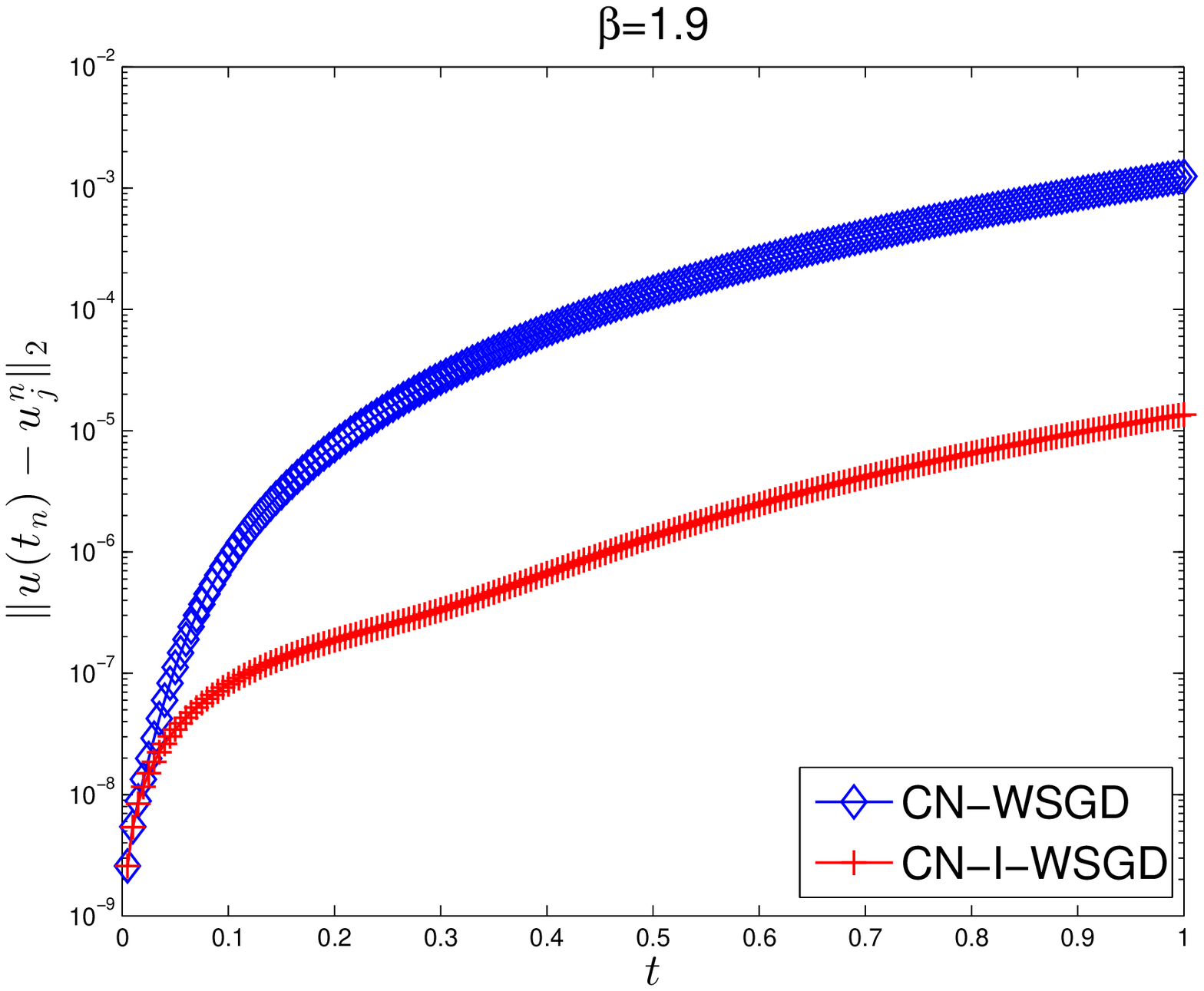}
	\caption{Comparison of asymptotic errors in $L_2$-norm for the CN-WSGD scheme \eqref{td-s} and the improved scheme CN-I-WSGD for the SpFDE \eqref{td} with left-sided fractional derivative (Example \ref{exm3}). The step-sizes are taken as $h=2^{-5}$ and $\tau=10^{-3}$.}
	\label{fig:soltd}
\end{figure}

Table \ref{table5} shows that the use of Algorithm \ref{algo1}  can greatly improve the accuracy and convergence rate of the fully discretized scheme CN-WSGD for solving the time-dependent problems. Compared to the low accuracy of numerical solutions produced by the scheme CN-WSGD \eqref{td-s}, the numerical solution from the improved scheme CN-I-WSGD  is of second-order convergence for $\beta=1.4$ and $1.8$. Fig \ref{fig:soltd} shows the asymptotic error in $L_2$-norm of numerical solutions. It can be observed that the improved scheme CN-I-WSGD performs far better than the original scheme CN-WSGD without applying the Algorithm \ref{algo1}.

\section{Conclusion}
					
We  proposed an improved algorithm for fractional boundary value problems with non-smooth solution by applying the extrapolation technique to the WSGD scheme \eqref{scheme-1}-\eqref{scheme-2} and the FCD scheme \eqref{scheme-c1}-\eqref{scheme-c2}. For some known structure of singularity, we  proved that the improved schemes I-WSGD and I-FCD can be of second-order convergence for non-smooth solution.  Numerical examples  show that the proposed algorithm \ref{algo1} and the improved schemes I-WSGD and I-FCD can significantly increase the accuracy and convergence rate of numerical solutions for fractional boundary value problems with one-sided fractional derivative or symmetric two-sided fractional derivatives with non-smooth solution; see Example \ref{exm1}-\ref{exm2}. Moreover, we  showed that the proposed algorithm can be successfully applied to the time-dependent problems with non-smooth solution and obtain highly accurate numerical solutions; see Example \ref{exm3}.

In the end, we give some remarks on the proposed algorithm. First, in this work, we  focused on dealing with the leading weak singularity for non-smooth solution.  To further increase the accuracy, especially for those problems with small fractional order $\beta$, the algorithm can be repeatedly applied for hierarchical singular terms. Second, though we have just applied the proposed algorithm to the WSGD scheme and the FCD scheme, it can be readily used to improve the accuracy of most finite difference schemes for SpFDEs.  Moreover, we have only studied the special cases for $\theta=0, 1/2, 1$ in this work. For the general case of $\theta\in(0,1),$ due to the difficulty of determining the function $\rho(\theta,\beta)$ in \eqref{eq:singularterm}, we do not consider the problem in this work.  In future work, we will extend  the proposed algorithm to  time-fractional differential equations and the general SpFDEs with non-smooth solution,  e.g.  variable-coefficient problems \cite{HaoPLC2016} and fractional advection diffusion equation \cite{Sousa2012}.

\section*{Acknowledgement}
We would
like to thank Dr. Sheng Chen  for helpful discussion  during the first author's  visiting in Purdue University     and thank Prof.  Zhi-Zhong Sun for proofreading the manuscript.
%\vskip 0.4cm
%{ \textbf{References}}

\end{document}